\documentclass[reqno,a4paper,11pt,oneside,final]{amsart}   
\usepackage{mathtools}

\usepackage{enumitem}
\setlist[enumerate,1]{label=(\roman*)}

\usepackage{xcolor}
\usepackage{booktabs}
\usepackage[english]{babel}
\usepackage{csquotes}
\usepackage[T1]{fontenc}
\usepackage{lmodern}

\usepackage[final,stretch=10]{microtype}
 
\usepackage[colorlinks, citecolor=hanblue,linkcolor=red]{hyperref}

\usepackage[
	style=authoryear-comp,
	sorting=nyt,
	url=true,
	doi=true,                
	isbn=false,              
	eprint=true,
	uniquename=init,         
	maxbibnames=99,          
	maxcitenames=3,
	uniquelist=minyear,
	dashed=false,            
	sortcites=false,          
	backend=biber,           
	safeinputenc,            
	date=year,               
]{biblatex}
\DeclareCiteCommand{\parencite}[\mkbibbrackets]{%
	\ifbibmacroundef{cite:init}{}{\usebibmacro{cite:init}}\usebibmacro{prenote}%
}{%
	\usebibmacro{citeindex}%
	\printtext[bibhyperref]{\usebibmacro{cite}}%
}{%
	\ifbibmacroundef{cite:init}{\multicitedelim}{}%
}{%
	\usebibmacro{postnote}%
}%

\let\cite\parencite

\DeclareFieldFormat{eprint:urn}{%
URN\addcolon\space
\ifhyperref
{\href{http://nbn-resolving.org/#1}{\nolinkurl{#1}}}
{\nolinkurl{#1}}}

\usepackage[noabbrev,nameinlink,capitalize]{cleveref}
\definecolor{hanblue}{rgb}{0.27, 0.42, 0.81}

\DeclarePairedDelimiter\parens()
\DeclarePairedDelimiter\bracks[]
\DeclarePairedDelimiter\braces\{\}

\DeclarePairedDelimiter{\abs}{\lvert}{\rvert}
\DeclarePairedDelimiter{\norm}{\lVert}{\rVert}

\DeclarePairedDelimiter{\seq}()

\providecommand\given{\nonscript\;\delimsize|\nonscript\;\mathopen{}}
\DeclarePairedDelimiterX\set[1]\{\}{#1}

\DeclarePairedDelimiterX\dual[2]{\langle}{\rangle}{#1,#2}
\DeclarePairedDelimiterX\ddual[2]{\llangle}{\rrangle}{#1,#2}
\DeclarePairedDelimiterX\scp[2](){#1,#2}

\newcommand\oo{o}

\DeclareMathOperator{\dom}{dom}
\DeclareMathOperator*{\argmin}{arg\,min}
\DeclareMathOperator{\sign}{sign}

\newcommand{\N}{\mathbb{N}}
\newcommand{\Z}{\mathbb{Z}}
\newcommand{\R}{\mathbb{R}}

\newcommand{\MM}{\mathcal{M}}

\newcommand{\U}{\mathcal{U}}
\newcommand{\Uc}{\mathcal{U}^*}

\renewcommand\d{\mathop{}\!\mathrm{d}}

\newcommand{\eps}{\varepsilon}

\newcommand{\weakstar}{\stackrel{*}{\rightharpoonup}}

\newcommand{\HH}{\mathcal{H}}

\newcommand{\KLfull}{Polyak--Łojasiewicz--Kurdyka}
\newcommand{\KL}{PŁK}

\newcommand\Uad{U_{\mathrm{ad}}}

\DeclareFontFamily{OMX}{MnSymbolE}{}
\DeclareSymbolFont{MnLargeSymbols}{OMX}{MnSymbolE}{m}{n}
\SetSymbolFont{MnLargeSymbols}{bold}{OMX}{MnSymbolE}{b}{n}
\DeclareFontShape{OMX}{MnSymbolE}{m}{n}{
    <-6>  MnSymbolE5
   <6-7>  MnSymbolE6
   <7-8>  MnSymbolE7
   <8-9>  MnSymbolE8
   <9-10> MnSymbolE9
  <10-12> MnSymbolE10
  <12->   MnSymbolE12
}{}
\DeclareFontShape{OMX}{MnSymbolE}{b}{n}{
    <-6>  MnSymbolE-Bold5
   <6-7>  MnSymbolE-Bold6
   <7-8>  MnSymbolE-Bold7
   <8-9>  MnSymbolE-Bold8
   <9-10> MnSymbolE-Bold9
  <10-12> MnSymbolE-Bold10
  <12->   MnSymbolE-Bold12
}{}
 
\makeatletter
\let\llangle\@undefined
\let\rrangle\@undefined
\makeatother
\DeclareMathDelimiter{\llangle}{\mathopen}%
                     {MnLargeSymbols}{'164}{MnLargeSymbols}{'164}
\DeclareMathDelimiter{\rrangle}{\mathclose}%
                     {MnLargeSymbols}{'171}{MnLargeSymbols}{'171}

\RequirePackage{calc}
\newcommand{\mrep}[3][l]{%
	\ifmmode%
	\mymathpalette\mrepinternal{#1}{#2}{#3}%
	\else%
	\makebox[\widthof{#3}][#1]{\vphantom{#3}#2}%
	\fi%
}

\def\mymathpalette#1#2#3#4{%
	\mathchoice
	{#1\displaystyle{#2}{#3}{#4}}%
	{#1\textstyle{#2}{#3}{#4}}%
	{#1\scriptstyle{#2}{#3}{#4}}%
	{#1\scriptscriptstyle{#2}{#3}{#4}}%
}

\def\mrepinternal#1#2#3#4{%
	\makebox[\widthof{$#1#4$}][#2]{$#1\vphantom{#4}{#3}$}%
}

\newcommand{\mscLink}[1]{\href{https://zbmath.org/classification/?q=cc:#1}{#1}}

\numberwithin{equation}{section}
 
\theoremstyle{plain}
\newtheorem{theorem}{Theorem}[section]
\newtheorem{lemma}[theorem]{Lemma}
\newtheorem{corollary}[theorem]{Corollary}
\newtheorem{proposition}[theorem]{Proposition}
\newtheorem{algorithm}[theorem]{Algorithm}
\theoremstyle{definition} 
\newtheorem{assumption}[theorem]{Assumption}
\newtheorem{definition}[theorem]{Definition}
\newtheorem{example}[theorem]{Example}
\newtheorem{remark}[theorem]{Remark}
 
\def\cleverreffix#1{\AddToHook{env/#1/begin}{\crefalias{theorem}{#1}}}
\cleverreffix{lemma}
\cleverreffix{corollary}
\cleverreffix{proposition}
\cleverreffix{algorithm}
\cleverreffix{assumption}
\cleverreffix{definition}
\cleverreffix{example}
\cleverreffix{remark}

\begin{document}
\title[]{Proximal gradient methods in Banach spaces}

\pagestyle{myheadings}

 \author[G. Wachsmuth]{Gerd Wachsmuth}
\author[D. Walter]{Daniel Walter}

 \address[Gerd Wachsmuth]{Institute of Mathematics, BTU Cottbus–Senf\-ten\-berg, 03046 Cottbus, Germany}

\address[Daniel Walter]{Institut f\"ur Mathematik, Humboldt-Universit\"at zu Berlin, 10117 Berlin, Germany}

\email[Gerd Wachsmuth]{gerd.wachsmuth@b-tu.de}
 \email[Daniel Walter]{daniel.walter@hu-berlin.de}

\begin{abstract}
\small{ %
Proximal gradient methods are a popular tool for the solution of structured, nonsmooth minimization problems. In this work, we investigate an extension of the former to general Banach spaces and provide worst-case convergence rates for, both, convex and nonconvex, problem instances. Moreover, assuming additional regularity properties of stationary points, linear rates of convergence are derived. The theoretical results are illustrated for bang-bang type optimal control problems with partial differential equations which we study in the space of Radon measures. An efficient implementation of the resulting $L^1$-proximal gradient method is given and its performance is compared to standard $L^2$-proximal gradient as well as Frank-Wolfe methods. The paper is complemented by discussing the relationship among different regularity properties as well as by providing a novel characterization of the \KLfull\ property via second-order conditions involving weak* second subderivatives.
 
\vskip .3truecm \noindent Keywords:
Optimization in Banach spaces,
\KLfull\ property,
second-order optimality conditions,
second subderivatives,
strong metric subregularity.
 
  \vskip.1truecm \noindent 2020 Mathematics Subject Classification:
  \mscLink{46E27}, \mscLink{49K30}, \mscLink{49J52}, \mscLink{49J53}
}
\end{abstract}
 
\maketitle
\section{Introduction}\label{sec:intro}
 In this paper, we consider composite minimization problems
\begin{equation} \label{def:problem}
    \min_{u \in \U} J(u) \coloneqq \left \lbrack f (u)+ g(u)   \right \rbrack, \tag{$\mathcal{P}$}
\end{equation}
where~$f$ is a smooth, but not necessarily convex, function, while $g$ denotes a convex, but potentially nonsmooth, regularizer. Problems of this form appear in a wide array of fields, ranging, e.g., from applications in data science to inverse problems and optimal control, with $f$ often modelling a fidelity term while $g$ encodes an a-priori known structural information on the sought-for solution.

If $\U$ is a Hilbert space with induced norm $\norm{\cdot}_{\U}=\sqrt{(\cdot,\cdot)_{\U}}$, the first-order necessary optimality condition of \eqref{def:problem} is given by      
\begin{equation*}
\bar{u}= \operatorname{Prox}_{g,L}\left(\bar{u}-\frac{1}{L} \nabla f(\bar{u})\right),    \quad \text{where} \quad \operatorname{Prox}_{g,L}(u)= \argmin_{v \in \U} \left\lbrack g(v)+ \frac{L}{2} \norm{v-u}^2_{\U} \right \rbrack
\end{equation*}
denotes the \textit{proximal operator} associated to $(1/L) g$ for some arbitrary $L>0$. Applying a fixed-point iteration yields the celebrated \textit{proximal gradient method}
\begin{equation} \label{eq:Hilbprox}
    u_{k+1} = \operatorname{Prox}_{g,L_k}\left(u_k-\frac{1}{L_k} \nabla f(u_k)\right)= \argmin_{v \in \U} \left \lbrack (\nabla f(u_k),v-u_k)_{\U}+ g(v)+ \frac{L_k}{2} \norm{v-u_k}^2_{\U} \right \rbrack
\end{equation}
for a sequence of suitable stepsizes $L_k>0$. Originating from the seminal works \cite{Moreau1965, Martinet1970, FukushimaMine1981}, proximal gradient algorithms have received tremendous attention over the last decades. This can be attributed to its simplicity, merely requiring a first-order oracle for $f$ as well as an efficient way to evaluate the proximal operator associated to $g$. For a variety of interesting and practically relevant examples, the latter is even possible in closed form, establishing \eqref{eq:Hilbprox} as an attractive tool for the approximate solution of large scale, structured minimization problems. Moreover, many modern optimization methods build upon the proximal framework as a basis, further underlining the significance of the latter. For a modern, non-exhaustive
exposition see, e.g., \cite{BauschkeCombettes2011,Parikh2014,Bubeck2015,Beck2017,Teboulle2018} as well as the references therein.
A pivotal role in the analysis of Prox-like methods is taken on by regularity conditions at stationary points $\bar{u}$ such as variants of the \textit{\KLfull\ property}, \cite{Polyak1963,Lojasiewicz1962,Lojasiewicz1965, Kurdyka1998,Bolte2007}, i.e.,
\begin{equation}
J(u)-J(\bar u)\leq
				\frac{\lambda}{2} \inf\set{\norm{\nabla f(u)+ \xi}^2_{\Uc} \given \xi \in \partial g(u)}.
\end{equation}
It enables the derivation of stronger convergence results and faster rates,
see, e.g., \cite{Karimi2016,Drusvyatskiy2018,Necoara2019,MehlitznonLipschitz,Bento2025}.

Within the present work, we investigate an analogue of \eqref{eq:Hilbprox}, merely assuming that $\U$ is isomorphic to the topological dual of a separable Banach space $Y$ as well as $\nabla f(u) \in Y$. More in detail, we consider
\begin{equation} \label{eq:proxintro}
    u_{k+1} \in \argmin_{v \in \U} \left \lbrack \dual{\nabla f(u_k)}{v-u_k}_{\U}+ g(v)+ \frac{L_k}{2} \norm{v-u_k}^2_{\U} \right \rbrack
\end{equation}
where $\dual{\cdot}{\cdot}_{\U}$ denotes the duality pairing between $Y$ and $\U$. While we clearly recover \eqref{eq:Hilbprox} via the canonical isomorphism if $\U$ is a Hilbert space, we show that the extension to Banach spaces is crucial, yielding new, implementable algorithms for challenging problems together with fast convergence guarantees. 
More in detail, our contributions are threefold:
\begin{enumerate}
    \item For general $f$ with Lipschitz-continuous gradient, we show that the proposed algorithm is well-defined and provide worst-case, sublinear rates of convergence for the minimal subgradient norm encountered
    \begin{equation*}
                 \inf_{\ell \in \{0, \dots,k \}}\inf\set{\norm{\nabla f(u_{\ell+1})+ \xi}_{\Uc} \given \xi \in \partial g(u_{\ell+1}) }
                 ,
    \end{equation*}
    which serves as a measure of stationarity,
    as well as, if $f$ is convex, for the objective functional residual
    \begin{equation*}
        r_J(u_k)=J(u_k)- \min_{u \in \U} J(u)
        ,
    \end{equation*}
    see \cref{thm:convexsublin,thm:convergencenonconv}, respectively.
If, moreover, $\seq{u_k}_k$ converges strongly towards an element $\bar{u}$ on a subsequence and if $\bar{u}$ satisfies additional regularity assumptions, linear convergence for, both, objective functional values as well as iterates follows, see \cref{thm:convergenceresults,thm:convexlinear}. More in detail, for convex $f$, we will rely on the \textit{Quadratic Growth} (QG) property
\begin{equation*}
        \frac{\theta}{2} \norm{u-\bar{u}}^2_\U \leq J(u)-J(\bar u) \quad \forall u \in B_R(\bar u) \text{ with } J(u) \le J(\bar u) + \eta
        ,
\end{equation*}
while we require the \textit{\KLfull} (\KL) condition
\begin{equation*}
    				J(u)-J(\bar u)
				\leq
				\frac{\lambda}{2} \inf\set{\norm{\nabla f(u)+ \xi}^2_X \given \xi \in \partial g(u) \cap X}
				\quad \forall u \in B_R(\bar u) \text{ with } J(u) \le J(\bar u) + \eta
\end{equation*}
with $X \in \{\U^*, Y\}$ in the general, nonconvex case.
    \item It is evident that QG and \KL~are not equivalent in general, see \Cref{ex:growth_and_KL,ex:growth_and_KL_2}, but the latter implies the former under rather weak conditions, see \cref{prop:KL_implies_growth,prop:KL_implies_growth_take_2}. Assuming additional regularity, it is well-known that QG is equivalent to the \textit{sufficient second-order} condition (SSC)
        \begin{equation*}
				f''(\bar u) h^2 + g''(\bar u, -\nabla f(\bar u); h)
				>
				0
				\qquad
				\forall h \in \U \setminus \set{0},
			\end{equation*}
involving the Hessian of $f$ as well as the second subderivative of $g$, together with a nondegeneracy condition. Similarly, equivalence between \KL~and \textit{strong metric subregularity} (SMS)
\begin{equation*}
\norm{ u - \bar u }_\U				
				\le
					\frac 1 \nu \inf\set{ \norm{ \nabla f(u) + \xi }_{X} \given \xi \in \partial g(u) \cap X }
				\quad
				\forall u \in B_R(\bar u) \text{ with } J(u) \le J(\bar u) + \eta.
\end{equation*}
can be established, see \Cref{lem:SMS_implies_KL,lem:KL_and_QG_imply_SMS}. Naturally, this raises question under which additional condition all four concepts are equivalent. Our main result in this regard, \cref{thm:second_order_SMS}, gives a, to the best of our knowledge, new sufficient regularity assumption on the nonsmooth part $g$ such that SSC implies SMS, thus closing the circle.
    \item As a particular example, we consider optimal control problems with partial differential equations (PDEs) 
    \begin{equation*}
      \min_{u \in L^2(\Omega)} \frac12 \norm{S(u) - y_d}_{L^2(\Omega)}^2 \quad \text{s.t.} \quad  u \in \Uad \coloneqq \set{ u \in L^2 (\Omega) \given u_a \leq u \leq u_b}  
  ,
\end{equation*}
where $\Omega \subset \R^d$, $d \in \set{1,2,3}$,
is a bounded Lipschitz domain and $u_a < u_b$ are constants.
The choice of the control $u$ influences the behavior of the state variable $y=S(u)$ via a semilinear, elliptic PDE
\begin{equation*}
	-\Delta y + a(y) = u \quad\text{in } \Omega,
	\qquad
	y = 0 \quad\text{on } \partial\Omega.
\end{equation*}
Problems of this form favor \textit{bang-bang}-type minimizers, i.e., functions $\bar{u} \in \Uad$ with $\bar{u}(x)=u_a$ or $\bar{u}(x)=u_b$ a.e. in $\Omega$, as the associated first-order necessary condition implies
\begin{equation*}
    \bar{u}(x) \in
    \begin{cases}
      \{u_a\}   & \text{if } \bar{p}(x)>0, \\
      \{u_b\}   & \text{if } \bar{p}(x)<0,  \\
      [u_a,u_b] & \text{if } \bar{p}(x)=0
    \end{cases}
    \quad \text{for a.e. }x \in \Omega,
\end{equation*}
where $\bar{p}=\nabla f(\bar{u})$ satisfies the adjoint equation
    \begin{equation*}
	-\Delta \bar p + a'(\bar y) p = \bar y - y_d \quad\text{in }\Omega,
	\qquad
	\bar p = 0 \quad\text{on } \partial\Omega.
\end{equation*}
Setting
\begin{equation*}
    f(u)=\frac12 \norm{S(u) - y_d}_{L^2(\Omega)}^2, \quad g(u)= \delta_{\Uad}(u)= \begin{cases}
        0 & \text{if }u \in \Uad, \\
        +\infty & \text{else},
    \end{cases}
\end{equation*}
we study the particular choice of $\U=\mathcal{M}(\Omega)$, the space of Radon measures on $\Omega$. Relying on the aforementioned abstract results, we prove equivalence of SSC, QG, \KL~  and SMS for $X=Y=C_0(\Omega)$ under the common regularity assumption
    \begin{equation} \label{eq:regbangintro}
      \lambda^d \left(\set{x \in \Omega \given |\bar{p}(x)|\leq \eps}\right) \leq C \eps \quad\forall \eps > 0
      \quad\text{and}\quad
      \nabla \bar{p}(x) \neq 0 \quad \text{on }~  \set{x \in \Omega \given |\bar{p}(x)|=0}
      ,
    \end{equation}
where $\lambda^d$ denotes the $d$-dimensional Lebesgue measure. Subsequently, we give an efficient implementation of the resulting $L^1$-proximal gradient method
    \begin{equation*}
    u_{k+1} \in \argmin_{v \in \Uad} \left \lbrack \int_\Omega p_k (v-u_k) \d \lambda^d(x)+ \frac{L_k}{2} \norm{v-u_k}^2_{L^1(\Omega)} \right \rbrack,
    \end{equation*}
where $p_k$ is again obtained by solving an adjoint equation. Its performance as well as the structure of the obtained iterates is compared to, both, the ``standard'' $L^2$ projected gradient method as well as a Frank-Wolfe method with backtracking step-size choice.
\end{enumerate}
In particular, our explicit, infinite-dimensional example highlights the necessity of an informed choice of the variable space $\U$ in order to exploit expectable and favourable geometric properties of $J$ in the derivation and analysis of solution algorithms.
\subsection{Related work} \label{subsec:related}
Proximal gradient methods are a vivid and fast-developing area of research. Hence, giving a full account of all related research goes beyond the scope of the current paper and we restrict ourselves to a tangible selection as well as the references contained within. Proximal gradient methods are a special case of a forward-backward splitting type algorithm applied to the necessary optimality condition
\begin{equation*}
    0 \in \nabla f(\bar{u}) + \partial g(\bar{u}).
\end{equation*}
 Abstract generalizations of the latter to Banach spaces are, e.g., discussed in \cite{Lopez2012,Cholamjiak2019}. In this regard, the work closest to ours is \cite{Bredies2008} which proposes a similar procedure to \eqref{eq:proxintro}. Subsequently, this idea was adopted by \cite{Guan2015,Guan2021} where the discrepancy term is replaced by the Bregman distance associated to the respective power of the norm. We refer also to \cite{Lazzaretti2022} for an application to problems in variable exponent Lebesgue spaces. In comparison to the present work, all of these require stronger assumptions on $\U$, in particular reflexivity. We emphasize that the present work does not rely on this restriction which would, in particular, exclude the bang-bang example which is studied in the nonreflexive space $\U= \mathcal{M}(\Omega)$. Moreover, we point out \cite{Halabi2017} in which the authors study \eqref{eq:proxintro} for finite-dimensional spaces $\U$ and non-euclidean norms, noting that a proper choice of the norm entails benefits on the convergence behaviour. Similar observations have been made prior for specific instances, see e.g. \cite{Boyer2014,Daspremont2018,Kelner2014}.

 As already stated earlier, interest in regularity concepts beyond strong convexity, such as QG, \KL, SMS and SSC, has increased tremendously over the last years. Clearly, depending on the application, particular concepts are more desirable than others since they, e.g., are easier to verify or lead to elegant convergence proofs. As a consequence, there have been several works exploring the relationships among them. We point out, e.g., \cite{CorellaLe2024,BorchardWachsmuth2024,LiaoDingZheng2024,LiMengYang2023,MohammadiSarabi2020} as well as the references therein for a nonexhaustive overview.

Finally, regarding optimal control of PDEs with bang-bang controls, the first work concerning regularity conditions in the form of second-order conditions is \cite{Casas2012},
in which quadratic growth w.r.t.\ the state variable $y$ was shown.
In order to obtain growth for the control, regularity assumptions akin to \eqref{eq:regbangintro} were used in
\cite{CasasWachsmuthWachsmuth2017:2}.
Afterwards, the analysis was refined and no-gap second-order conditions were provided,
see \cite{ChristofWachsmuth2017:1,wachsmuth2}.
Finally, we mention the contributions
\cite{Wachsmuth2008,DeckelnickHinze2012,WachsmuthWachsmuth2009,vonDaniels2016},
which also
used conditions similar to \eqref{eq:regbangintro} to analyze problems
with bang-bang controls.

\subsection{Outline} \label{subsec:outline}
Necessary notation and concepts are introduced in \Cref{sec:notation}. \Cref{sec:problen} formally introduces the proposed method and proves its wellposedness while worst-case convergence rates for both convex and nonconvex problems are derived in \Cref{sec:worstcase}. Linear convergence rates in different settings and for different regularity conditions are presented in \Cref{sec:improved}. \Cref{sec:KL} explores the relationship between QG, \KL, SMS and SSC. Finally, we illustrate the abstract results on the particular example of bang-bang control problems in \Cref{sec:bangbang} and present numeric experiments.   

\section{Notation \& standing assumptions} \label{sec:notation}
We define
$\R^+ \coloneqq [0,\infty)$, $\bar\R^+ \coloneqq [0,\infty]$, and $\bar\R := (-\infty,\infty]$.
Further, $\Z$ is the set of integers and $\N$ is the set of non-negative integers.
By $\lambda^d$, we refer to the Lebesgue measure on $\R^d$.

The following properties are tacitly assumed throughout the paper.
\begin{assumption} \label{ass:setup}
Assume that:
\begin{enumerate}[label=\textbf{A\arabic*}]
  \item\label{ass:setup:1}
    We have $\U \cong Y^*$ where~$Y$ is a separable Banach space with norm $\norm{\cdot}_Y$. The associated duality pairing is denoted by~$\dual{\cdot}{\cdot}$ and we equip $\U$ with the canonical dual norm
    \begin{align*}
      \norm{u}_\U \coloneqq \sup_{\norm{y}_Y\leq 1}\dual{y}{u} \quad \forall u \in \U.
    \end{align*}
  \item\label{ass:setup:2}
    The function $g \colon \U \to \Bar{\R}^+ $ is proper, convex, and sequentially weak* lower semi-continuous.
  \item\label{ass:setup:3}
  The function $f \colon  \U \to \R^+$ is sequentially weak* lower semi-continuous on $\dom g$ and there is a neighborhood $\mathcal{N}_{\U}$ of $\dom g$ in $ \U$ as well as a mapping $\nabla f \colon \dom g \to Y $ such that $f$ is Gâteaux differentiable on $\mathcal{N}_\U$ with
    \begin{align*}
      f'(u)( \delta u) = \dual{ \nabla f(u)}{\delta u} \quad \forall u \in \dom g, \delta u \in \U.
    \end{align*}
    Moreover, $\nabla f$ is Lipschitz-continuous on $\dom g$, i.e., there is $L_{\nabla f} \ge 0$ such that
    \begin{align*}
      \norm{\nabla f(u_1)-\nabla f(u_2)}_Y \leq L_{\nabla f} \norm{u_1-u_2}_{\U} \quad \forall u_1, u_2 \in \dom g.
    \end{align*}
  \item\label{ass:setup:4}
    The objective functional is radially unbounded, i.e., for sequences~$\seq{u_k}_k \subset \U$ we have
    \begin{align*}
      \norm{u_k}_\U \rightarrow \infty \quad\Rightarrow\quad J(u_k) \rightarrow \infty.
    \end{align*}
\end{enumerate}
\end{assumption}
We further denote the duality pairing between $\U$ and its dual space $\Uc$ by $\ddual{\cdot}{\cdot}$ and define the canonical norm on $\Uc$ by
\begin{align*}
        \norm{\zeta}_{\Uc} \coloneqq \sup_{\norm{u}_\U\leq 1}\ddual{\zeta}{u} \quad \forall \zeta \in \Uc.
\end{align*}
Note that we have
\begin{equation*}
    \dual{\varphi}{u}=\ddual{\varphi}{u}, \quad \norm{y}_Y = \norm{y}_{\Uc} \quad \forall \varphi \in Y,~u \in \U 
\end{equation*}
due to the canonical and isometric embedding $Y \hookrightarrow \Uc=Y^{**} $. For a proper and convex function $h \colon \U \to \bar{\R}$, we define its subdifferential at $u \in \U$ by
\begin{equation*}
    \partial h(u) \coloneqq \set{\xi \in \Uc \given \ddual{\xi}{v-u}+h(u) \leq h(v) \quad \forall v \in \U}.
\end{equation*}
Note that $\partial h(u) \subset \Uc$ by definition.
At some places, we utilize $\partial h(u) \cap Y$ (again using the embedding $Y \hookrightarrow \Uc$)
for which we have the representation
\begin{equation*}
    \partial h(u) \cap Y = \set{\xi \in Y \given \dual{\xi}{v-u}+h(u) \leq h(v) \quad \forall v \in \U}.
\end{equation*}

Given $R>0$ as well as $\bar{u} \in \U$, let $B_R(\bar{u}) \coloneqq \set{u \in \U \given \norm{u-\Bar{u}}_{\U} \leq R }$ denote the closed ball of radius $R$ centered at $\bar{u}$. An element $\bar{u} \in \U$ is called a \text{local minimizer} of \eqref{def:problem} if we have
\begin{equation*}
    J(u)-J(\bar{u}) \geq 0 \quad \forall u \in B_R(\bar{u})
\end{equation*}
for some $R>0$. Similarly, we call $\bar{u}$ a \text{stationary point} of \eqref{def:problem} if $-\nabla f(\bar{u}) \in \partial g(\bar{u})$. Associated to $J$, we further define the residual
\begin{equation*}
    r_J(u)\coloneqq J(u)- \min \eqref{def:problem} \quad \forall u \in \U.
\end{equation*}
Throughout the paper, $u_0 \in \dom g$ will denote an arbitrary but fixed initial iterate.
By assumption \ref{ass:setup:4}, the associated sublevel set
\begin{equation*}
    E(u_0) \coloneqq \set{u \in \U\given J(u) \leq J(u_0)}
\end{equation*}
is bounded, i.e., there exists $M_{0}>0$ such that
\begin{equation*}
    \sup \set{\norm{u}_\U \given u \in E(u_0) } \leq M_{0}.
\end{equation*}
Finally, we introduce the associated dual gap functional as
\begin{equation*}
    \Psi(u)= \max_{\norm{v}_{\U} \leq M_0} \left \lbrack \dual{\nabla f(u)}{u-v} + g(u)-g(v) \right \rbrack \quad \text{for all} \quad u \in E(u_0).
\end{equation*}
If $f$ is convex, we have $r_J(u)\leq \Psi(u) $ for all $u \in E(u_0)$. 

Following \cite{ChristofWachsmuth2017:1,wachsmuth2,BorchardWachsmuth2024},
we introduce the weak* second subderivative.
\begin{definition}
	\label{def:weak_star_subderivative}
	Let $u \in \dom g$ and $\xi \in Y$ be given.
	Then the weak* second subderivative
	$g''(u, \xi; \cdot) \colon \U \to [-\infty, \infty]$
	of $g$ at $u$ for $\xi$ is defined by
	\begin{equation*}
		g''(u, \xi; h)
		:=
		\inf
		\set*{
			\liminf_{k \to \infty} \frac{g(u + t_k h_k) - g(u) - t_k \dual{\xi}{h_k}}{t_k^2/2}
			\given
			t_k \searrow 0,
			h_k \weakstar h
		}
		.
	\end{equation*}
\end{definition}
For more information concerning the weak* second subderivative,
we refer to the discussions in
\cite{ChristofWachsmuth2017:1,wachsmuth2,BorchardWachsmuth2024}.
On $\U = \R^n$, the (weak*) second subderivative
is a classical object, see \cite[13.3~Definition]{RockafellarWets1998}.
At this point, we just mention that,
under some further, mild assumptions on $f$ and $g$,
one can prove that the second-order condition
\begin{equation*}
	f''(u) h^2 + g''(u, \xi; h^2) > 0
	\qquad\forall h \in \U
\end{equation*}
for $u \in \dom g$ and $\xi \in Y$
is equivalent to the quadratic growth condition
\begin{equation*}
	J(u) \ge J(\bar u) + \frac{\theta}{2} \norm{u - \bar u}_{\U}^2
	\qquad
	\forall u \in B_R(\bar u)
\end{equation*}
for some constants $\theta, R > 0$,
see
\cite[Theorem~4.5]{ChristofWachsmuth2017:1},
\cite[Theorem~2.10]{wachsmuth2},
or
\cite[Theorem~2.20]{BorchardWachsmuth2024}.

 \section{A proximal gradient method} \label{sec:problen}
 By standard arguments, we obtain the existence of minimizers to Problem \eqref{def:problem} as well as first-order necessary optimality conditions.  
\begin{theorem} \label{thm:existence}
There exists at least one global minimizer of Problem \eqref{def:problem}. Moreover, the set of minimizers to Problem \eqref{def:problem} is bounded in $\U$. Every local minimizer $\bar{u}$ of Problem \eqref{def:problem} satisfies $-\nabla f(\bar{u}) \in \partial g(\bar{u})$.
\end{theorem}
Note that $\bar{u}$ is a stationary point of \eqref{def:problem} if and only if
\begin{equation} \label{def:stationaritymeasures}
    \inf\set{\norm{\nabla f(\bar{u})+ \xi}_{\Uc} \given \partial g(\bar{u}) }=0,
\end{equation}
which will serve as a quantitative measure of stationarity throughout the paper.

Given~$L>0$, we define the model function
$\mathcal{Q}_L \colon \U \times \dom g \to \bar\R$
via
\begin{align*}
    \mathcal{Q}_L (v,u)
    \coloneqq
    f(u)+ \dual{\nabla f(u)}{v-u}+ g(v)+\frac{L}{2} \norm{v-u}_\U^2,
\end{align*} 
which can be interpreted as a partial linearization of $J$ around $u$ with a quadratic error term. The following two propositions address the minimization of this auxiliary functional w.r.t.\ $v$.
\begin{proposition} \label{prop:wellposedalg}
Let~$L>0$ and $u \in \dom g$ be given. Then there exists a minimizer $\bar{v}$ of
\begin{align} \label{def:proxprobaux}
    \inf_{v\in \U} \mathcal{Q}_L (v,u).
\end{align}
Moreover, the set of minimizers $\mathcal{S}$ is convex and we have
\begin{equation*}
\norm{v_1-u}_\U=\norm{v_2-u}_\U, \quad g((1-\sigma)v_1+\sigma v_2)=(1-\sigma)g(v_1)+\sigma g(v_2)
\end{equation*}
for all $v_1,v_2 \in \mathcal{S}$ and $\sigma \in [0,1]$.
\end{proposition}
\begin{proof}
Let~$\seq{v_j}_j \subset \U$ denote an infimizing sequence for~\eqref{def:proxprobaux}, i.e., there holds
\begin{align*}
\mathcal{Q}_L (v_j,u) \rightarrow  \inf_{v\in \U} \mathcal{Q}_L (v,u) < +\infty 
.
\end{align*}
For $j$ large enough, we then have
\begin{align*} \label{eq:existenceproxaux1}
  J(u) + 1
  =
  \mathcal{Q}_L (u,u) + 1
  &\geq \mathcal{Q}_L (v_j,u)
  =
  f(u)+ \dual{\nabla f(u)}{v_j-u}+ g(v_j)+\frac{L}{2} \norm{v_j-u}_\U^2 \\
  & \geq f(u) + g(v_j) + \left(-\norm{\nabla f(u)}_Y+\frac{L}{2}\norm{v_j-u}_\U  \right)\norm{v_j-u}_\U
\end{align*}
as well as
\begin{equation*}
    \norm{v_j}_\U \leq \norm{v_j-u}_\U+ \norm{u}_\U.
\end{equation*}
Since the convex function $g$ possesses a bounded and affine minorant,
$\seq{v_j}_j$ is bounded in $\U$ and thus admits a weakly*-convergent subsequence, denoted by the same subscript, with limit $\bar{v}$. The proof is finished noting that $\mathcal{Q}_L (\cdot,u)$ is sequentially weak*-lower semicontinuous, i.e.,
 \begin{align*}
   \inf_{v\in \U} \mathcal{Q}_L (v,u) \leq  \mathcal{Q}_L (\bar{v},u) \leq   \liminf_{j \rightarrow \infty} \mathcal{Q}_L (v_j,u)=\inf_{v\in \U} \mathcal{Q}_L (v,u).
 \end{align*}
 As a consequence, $\bar{v}$ is a minimizer of \eqref{def:proxprobaux}. The remaining statements follow by convexity of $g$ as well as strict convexity of the square.
\end{proof}
We emphasize that the mapping $v \mapsto \norm{v-u}^2_\U $ is, in general, not strictly convex, i.e., the set of minimizers in \eqref{def:proxprobaux} is typically not a singleton.
\begin{proposition} \label{prop:firstorderpox}
    Let $L > 0$ and $u \in \dom g$ be given.
An element $\bar{v} \in \U$ is a minimizer of \eqref{def:proxprobaux} if and only if there exists $\xi \in \Uc$ with
\begin{align*}
       -\nabla f(u)- \xi \in \frac{L}{2} \partial \norm{\cdot-u}_\U^2 (\bar{v}) \quad \text{and} \quad \xi \in \partial g(\bar{v}).
       \end{align*}
\end{proposition}
\begin{proof}
    Since \eqref{def:proxprobaux} is a convex problem, $\bar{v}$ is a minimizer of \eqref{def:proxprobaux} if and only if
    \begin{align*}
        0
        \in
        \partial \mathcal{Q}_L (\cdot,u)(\bar{v})
        &= \{\nabla f(u)\}+ \partial \left( g(\cdot)+ \frac{L}{2}\norm{\cdot-u}_\U^2 \right)(\bar{v})
        \\
        &=\{\nabla f(u)\}+ \partial g(\bar{v}) + \frac{L}{2} \partial \norm{\cdot-u}_\U^2 (\bar{v}),
    \end{align*}
    where we apply the subdifferential sum rule twice noting that
    the linear functional $\nabla f(\bar u)$ and
    the norm are continuous.
    This finishes the proof. 
\end{proof}
\begin{remark}
As we will see in the following, it would be favourable to have $\xi \in Y \hookrightarrow \Uc $. However, this additional regularity cannot be guaranteed in general.
This would be true, e.g.,
if we are able to apply the subdifferential sum rule
in the weak* topology of $\U$,
but this requires the very strong assumption that $g$
possesses a point of continuity
w.r.t.\ the weak* topology of $\U$. For, both, positive and negative examples in this regard,  we refer to \cref{lem:subgradients_predual_bangbang,ex:continuous_subgradients} below.
   
\end{remark}
The following descent lemma gives motivation to the auxiliary problem \eqref{def:proxprobaux}. Its proof follows a standard Taylor approximation argument and is omitted for the sake of brevity.
 \begin{proposition} \label{prop:descent}
    For all~$u_1, u_2 \in \dom g $ and~$L\geq L_{\nabla f}$ we have
    \begin{equation*}
        f(u_1)-f(u_2) \leq \dual{\nabla f(u_2)}{u_1-u_2}+ \frac{L}{2} \norm{u_1-u_2}^2_\U 
    \end{equation*}
    and thus
        \begin{equation} \label{eq:descentineq}
    J(u_1) \leq \mathcal{Q}_L(u_1,u_2).
\end{equation}
\end{proposition}
\Cref{prop:descent} motivates the following iterative procedure.
\begin{algorithm}
  \label{alg:prox}
  We fix a starting point $u_0 \in \dom g$
  and constants $0 < L_a \le L_b < \infty$.
  For each $k \in \N$, we select
  \begin{align} \label{eq:prox}
    u_{k+1} \in \argmin_{v \in \U } \mathcal{Q}_{L_k} (v,u_k)= \argmin_{v \in \U } \left\lbrack \dual{\nabla f(u_k)}{v}+ g(v)+ \frac{L_k}{2} \norm{v-u_k}^2_\U   \right \rbrack,
  \end{align}
  where $\seq{L_k}_k$ is a sequence of stepsize parameters which are chosen such that we have
  \begin{align}\label{eq:condonL}
    L_k \in [L_a, L_b]
    \quad\text{and}\quad
    J(u_{k+1})\leq \mathcal{Q}_{L_k} (u_{k+1},u_k)
    .
  \end{align}
\end{algorithm}
Note that we do not require $L_k \geq L_{\nabla f}$,
but merely assume that the analogue of \eqref{eq:descentineq} is satisfied at every iteration.
Due to \cref{prop:descent}, \Cref{alg:prox} is well-defined
if $L_b \ge L_{\nabla f}$.
\begin{remark} \label{rem:stepsizecho}
Since $L_{\nabla f}$ is often unknown in practice and in order to better exploit local properties of $\nabla f$, $L_k$ can, e.g., be chosen by a backtracking procedure.
More in detail,
given scaling factors $\gamma_l >1$ and $\gamma_s \in (0,1)$ as well as an initial and a minimal stepsize $L_{0,0}, L_a>0$, respectively,
we choose the smallest $n_k \in \N$, such that
$L_k \coloneqq L_{k,0} \gamma^{n_k}_l $
and
$u_{k+1}$ from \eqref{eq:prox}
satisfy
\begin{equation*}
    J(u_{k+1}) \leq \mathcal{Q}_{L_k} (u_{k+1},u_k) .
\end{equation*}
Subsequently, we set $L_{k+1,0} \coloneqq \max\{L_k \gamma_s, L_a\}$ and proceed to the next iteration. Recalling \cref{prop:descent}, this procedure is well posed and we have $L_a \leq L_k \leq L_b \coloneqq L_{\nabla f} \gamma_l $.
\end{remark}
By construction, there holds
\begin{align*}
    J(u_{k+1})\leq \mathcal{Q}_L (u_{k+1},u_k) \leq  \mathcal{Q}_L (u_{k},u_k)=J(u_k).
\end{align*}
Hence, the sequence of objective values $\seq{J(u_k)}_k$ is decreasing.
In the remainder of this section, we collect some pertinent results for the algorithm described in \cref{alg:prox}. The first result addresses the per-iteration descent achieved by \cref{alg:prox}. 
\begin{lemma} \label{lem:descentestimate}
Let~$\seq{u_k}_k$ be generated by~\Cref{alg:prox}.
For all~$k \in \N$, there is~$\xi_k \in \partial g(u_{k+1}) $ such that we have
\begin{equation}
    \label{eq:firstorderaux}
    -\nabla f(u_k)- \xi_k \in \frac{L}{2} \partial \norm{\cdot-u_k}_\U^2 (u_{k+1})
    .
\end{equation}
Every such $\xi_k$ satisfies
\begin{subequations}
    \label{eq:stuff_for_xi}
    \begin{align}
        \label{eq:stuff_for_xi:1}
        L_k \norm{u_{k+1}-u_k}^2_{\U}
        &=
        \frac{1}{L_k} \norm{\nabla f(u_k)+ \xi_k}^2_{\Uc}
        = \ddual{\nabla f(u_k)+ \xi_k}{u_k-u_{k+1}}
        ,
        \\
        \label{eq:stuff_for_xi:2}
        J(u_k)-J(u_{k+1})
        &
        \geq \frac{L_k}{2} \norm{u_{k+1}-u_k}^2_{\U} 
        .
    \end{align}
\end{subequations}
\end{lemma}
\begin{proof}
    The existence of $\xi_k \in \partial g(u_{k+1})$
    satisfying
    \eqref{eq:firstorderaux}
  follows from \cref{prop:firstorderpox}.
  Due to the Fenchel-Young equality, this implies
  \begin{align*}
      0
      &=
      \frac{L_k}{2} \norm{u_{k+1}-u_k}^2_{\U}+ \frac{1}{2L_k} \norm{\nabla f(u_k)+ \xi_k}^2_{\Uc} - \ddual{\nabla f(u_k)+ \xi_k}{u_k-u_{k+1}}
      \\
      &\ge
      \frac{L_k}{2} \norm{u_{k+1}-u_k}_{\U}^2
      +
      \frac{1}{2 L_k }\norm{\nabla f(u_k)+ \xi_k}_{\Uc}^2
      -
      \norm{u_{k+1}-u_k}_{\U} \norm{\nabla f(u_k)+ \xi_k}_{\Uc}
      \\&=
      \frac{L_k}{2} \parens*{
          \norm{u_{k+1}-u_k}_{\U}
          -
          \frac{1}{L_k }\norm{\nabla f(u_k)+ \xi_k}_{\Uc}
      }^2
      .
  \end{align*}
  It is easy to check that this yields \eqref{eq:stuff_for_xi:1}.
  Using $J(u_{k+1}) \leq \mathcal{Q}_{L_k} (u_{k+1},u_k)$ due to \eqref{eq:condonL}, we have
  \begin{align*}
    J(u_k)-J(u_{k+1})
    &\geq \dual{\nabla f(u_k)}{u_k-u_{k+1}}+g(u_k)-g(u_{k+1})- \frac{L_k}{2} \norm{u_{k+1}-u_k}^2_\U \\
    &\geq \ddual{\nabla f(u_k)+ \xi_k}{u_k-u_{k+1}} - \frac{L_k}{2} \norm{u_{k+1}-u_k}^2_\U
    = \frac{L_k}{2} \norm{u_{k+1}-u_k}^2_\U,
  \end{align*}
  where the second inequality follows from $\xi_k \in \partial g(u_{k+1})$ and the equality uses
  \eqref{eq:stuff_for_xi:1}.
\end{proof}
Since $J$ is bounded from below,
the usual argument shows that
$\norm{u_{k+1}-u_k}_\U$ goes to zero.
\begin{corollary}\label{cor:convofsuccessive}
    Let~$\seq{u_k}_k$ be generated by~\Cref{alg:prox}. Then there holds~$\norm{u_{k+1}-u_k}_\U \rightarrow 0$ as~$k \rightarrow \infty$. 
\end{corollary} 
 \begin{proof}
    From \cref{lem:descentestimate} we get
    \begin{align*}
        \infty
        >
        r_J(u_0)
        \ge
        J(u_0) - J(u_{k})
        =
        \sum_{j = 0}^{k-1} J(u_j) - J(u_{j+1})
        \ge
        \frac{L_a}{2} \sum_{j = 0}^{k-1} \norm{u_{j+1} - u_j}_{\U}^2
        .
    \end{align*}
    Passing  to the limit $k \to \infty$ yields the claim.
\end{proof}

The following lemma ensures that the same holds true for the stationarity measure introduced in \eqref{def:stationaritymeasures}.
 \begin{lemma} \label{lem:estfordist}
 Let~$\seq{u_k}_k$ be generated by~\Cref{alg:prox}. Then there holds  
     \begin{equation*}
         \inf\set{\norm{\nabla f(u_{k+1})+ \xi}_{\Uc} \given \xi \in \partial g(u_{k+1}) }
         \leq
         ( L_{\nabla f}+ L_k)\norm{u_{k+1}-u_k}_\U
     \end{equation*}
     for all $k \in \N$,
     where~$L_{\nabla f}>0$ denotes the Lipschitz constant of~$\nabla f$.
 \end{lemma}
 \begin{proof}[Proof of \cref{lem:estfordist}]
    With $\xi_k$ from \cref{lem:descentestimate} we have
    \begin{align*}
     \inf\set{\norm{\nabla f(u_{k+1})+ \xi}_{\Uc} \given \xi \in \partial g(u_{k+1}) }
     & \leq \norm{\nabla f(u_{k+1})+ \xi_{k}}_{\Uc} \\
     & \leq \norm{\nabla f(u_{k})+ \nabla f(u_{k+1})}_{Y}+ \norm{\nabla f(u_{k})+ \xi_{k}}_{\Uc} \\
     & \leq ( L_{\nabla f}+ L_k) \norm{u_{k+1}-u_k},
    \end{align*}
    where we used the Lipschitz continuity of~$\nabla f$ and \eqref{eq:stuff_for_xi:1} in the final inequality.
 \end{proof}
\section{Worst-case convergence results} \label{sec:worstcase} 
This section is dedicated to the derivation of worst-case convergence results based on \cref{ass:setup}. For this purpose, and for the rest of the paper, let $u_0 \in \dom g$ denote a fixed initial iterate.
We start by showing that $\seq{u_k}_k$ admits weak* accumulation points and discuss the convergence of the associated functional values.
\begin{lemma} \label{lem:usefulstuff}
Let~$\seq{u_k}_k$ be generated by~\Cref{alg:prox}. Then~$\seq{u_k}_k$ admits at least one weak* accumulation point~$\Bar{u} \in \U$ and we have
\begin{align}  \label{eq:lowerbound}
    J(\bar{u}) \leq J(u_k) \quad \text{for all} \quad k \in \N.
\end{align}
If there is a subsequence~$\seq{u_{k_j}}_j$ with~$u_{k_j} \rightarrow \bar{u}$ strongly in~$\U$, we have~$\lim_{k \rightarrow \infty} J(u_k)=J(\bar{u})$.
 \end{lemma}
 \begin{proof}
 Recall that $\seq{J(u_k)}_k$ is monotonically decreasing and bounded from below. Hence, it converges towards a limit $\bar J$. The claimed existence of a weakly* convergent subsequence~$\seq{u_{k_j}}_j$ with limit~$\bar{u} \in \dom g$ follows immediately since~$\seq{u_k}_k$ is bounded in~$\U$. Due to the weak* lower semicontinuity of~$J$, we further have
\begin{align*}
    J(\bar{u}) \leq \liminf_{j \rightarrow \infty} J(u_{k_j})= \bar{J}
\end{align*}
and, thus,~\eqref{eq:lowerbound} follows due to the monotonicity of~$\seq{J(u_k)}_k$. It remains to show that we have
\begin{equation} \label{eq:auxlsc}
\lim_{k \rightarrow \infty} J(u_k)=J(\bar{u})  
\end{equation}
For this purpose, recall that~$f$ is continuous on $\dom g$, i.e., we have~$f(u_{k_j}) \rightarrow f(\bar{u})$, and~$u_{k+1}$ satisfies~\eqref{eq:prox}. The latter implies
\begin{align*}
    \dual{\nabla f(u_{k_j})}{u_{k_j+1}-u_{k_j}}&+ g(u_{k_j+1})+ \frac{L_{k_j}}{2} \norm{u_{k_j+1}-u_{k_j}}^2_\U \\ & \leq \dual{\nabla f(u_{k_j})}{\bar{u}-u_{k_j}}+ g(\bar{u})+ \frac{L_{k_j}}{2} \norm{\bar{u}-u_{k_j}}^2_\U.
\end{align*}
Taking the limit superior on both sides yields
\begin{align*}
    \limsup_{j \rightarrow \infty} g(u_{k_j}) \leq g(\bar{u})  \quad \text{and thus} \quad  \bar{J}=\limsup_{j \rightarrow \infty} J(u_{k_j}) \leq J(\bar{u})
\end{align*}
due to Lipschitz continuity of $\nabla f$, strong convergence of $\seq{u_{k_j}}_j$ as well as due to \cref{cor:convofsuccessive}. Combining both estimates, yields
\begin{equation*}
    \bar{J} =\lim_{k \to \infty} J(u_k)=\lim_{j \to \infty} J(u_{k_j})= J(\bar{u})
\end{equation*}
finishing the proof.
\end{proof}
Under an additional assumption, we can prove stationarity of the accumulation points.
\begin{lemma}
    \label{lem:stationarity_of_accumulation_points}
    Let~$\seq{u_k}_k$ be generated by~\Cref{alg:prox}.
    If
    \begin{enumerate}
        \item
            $u_{k_j} \to \bar u$ or
        \item
            $u_{k_j} \weakstar \bar u$ and $\nabla f \colon \dom g \to Y$ is sequentially weak*-to-strong continuous
    \end{enumerate}
    then $\bar u$ is stationary, i.e., $-\nabla f(\bar u) \in \partial g(\bar u)$.
\end{lemma}
\begin{proof}
    The assumptions guarantee
    $-\nabla f(u_{k_j}) \to -\nabla f(\bar u)$ in $Y$.
    Together with $u_{k_j} \weakstar \bar u$
    we can pass to the limit in $-\nabla f(u_{k_j}) \in \partial g(u_{k_j})$
    and obtain stationarity of $\bar u$.
\end{proof}
Next, we provide a rate of convergence for the stationarity measure.  
\begin{theorem} \label{thm:convergencenonconv}
    Let~$\seq{u_k}_k$ be generated by~\Cref{alg:prox}.
    Then there holds
           \begin{align*}
         \inf_{\ell \in \{0, \dots,k \}}\inf\set{\norm{\nabla f(u_{\ell+1})+ \xi}_{\Uc} \given \xi \in \partial g(u_{\ell+1}) } \leq \sqrt{\frac{ r_J(u_0)} {\nu(k+1)}} 
     \end{align*}
     for all $k \geq 0$ where
     \begin{equation*}
         \nu=\min\set*{
                \frac{L_a}{2( L_{\nabla f}+ L_a)^2},
                \frac{L_b}{2( L_{\nabla f}+ L_b)^2}}
                .
     \end{equation*}
\end{theorem}
\begin{proof}
    Combining \cref{lem:descentestimate,lem:estfordist}, we arrive at
    \begin{equation*}
        J(u_\ell) - J(u_{\ell+1})
        \geq \frac{L_\ell}{2} \norm{u_{\ell+1}-u_\ell}^2_{\U}
        \geq  \frac{L_\ell}{2( L_{\nabla f}+ L_\ell)^2} \inf\set{\norm{\nabla f(u_{\ell+1})+ \xi}_{\Uc} \given \xi \in \partial g(u_{\ell+1}) }^2
    \end{equation*}
    for every~$\ell \in \N$. Next, we sum both sides over~$\ell=0, \dots, k$ and note that
    \begin{equation*}
        \frac{L_\ell}{2( L_{\nabla f}+ L_\ell)^2} \geq \min\set*{
                \frac{L_a}{2( L_{\nabla f}+ L_a)^2}
                ,
                \frac{L_b}{2( L_{\nabla f}+ L_b)^2}}=\nu 
    \end{equation*}
    follows since $t \mapsto t/(L_{\nabla f}+t) $ is increasing on $[L_a, L_{\nabla f}]$ and decreasing on $[L_{\nabla f}, \infty]$. We obtain
    \begin{align*}
        r_J(u_0)
        &\geq
        J(u_0)-J(u_{k+1})
        =
        \sum^k_{\ell=0} \lbrack J(u_\ell)-J(u_{\ell+1}) \rbrack \\
        &\geq
        \nu (k+1) \inf_{\ell \in \{0, \dots,k \}}\inf\set{\norm{\nabla f(u_{\ell+1})+ \xi}_{\Uc} \given \xi \in \partial g(u_{\ell+1}) }^2.
    \end{align*}   
Rearranging yields the desired claim.
\end{proof}
\begin{remark} \label{rem:replaceUc}
We emphasize that $\partial g(\cdot)$ can replaced by $\partial g(\cdot) \cap Y$ in the statements of \cref{lem:estfordist} and \cref{thm:convergencenonconv}
if we can guarantee the existence of $\xi_k \in \partial g(u_{k+1}) \cap Y$ satisfying \eqref{eq:firstorderaux} for all $k \ge 0$,
cf.\ \cref{lem:descentestimate}.
Note that this yields stronger results due to
\begin{equation*}
    \inf\set{\norm{\nabla f(u_{k+1})+ \xi}_{\Uc} \given \xi \in \partial g(u_{k+1}) \cap Y }
    \ge
    \inf\set{\norm{\nabla f(u_{k+1})+ \xi}_{\Uc} \given \xi \in \partial g(u_{k+1}) }
    .
\end{equation*}
\end{remark}
\begin{remark} \label{rem:convdual}
\Cref{thm:convergencenonconv} also implies a convergence rate for the dual gap $\Psi(u_k)$ which is often considered as a measure of stationarity for (generalized) conditional gradient methods and which is, in many cases, easier to evaluate. Indeed, we obtain
    \begin{equation*}
       \Psi(u_k)= \max_{\norm{v}_\U \leq M_0} \left\lbrack \dual{\nabla f(u_k)}{u_k-v}+ g(u_k)-g(v) \right \rbrack \leq \ddual{\nabla f(u_k)+\xi}{u_k-v_k} \leq 2M_0 \norm{\nabla f(u_k)+\xi}_{\Uc} 
    \end{equation*}
    for every $\xi \in \partial g(u_k)$. Infimizing w.r.t.\ the latter and applying \cref{thm:convergencenonconv}, we  conclude
    \begin{equation*}
        \inf_{\ell \in \set{0,\ldots,k}} \Psi(u_{\ell+1}) \leq 2M_0 \sqrt{\frac{ r_J(u_0)} {\nu(k+1)}}. 
    \end{equation*}
\end{remark}
Next, we prove a sublinear rate of convergence for the residual as well as subsequential convergence towards minimizers if $f$ is convex. For this purpose, we compare the per-iteration descent achieved by \cref{alg:prox} with an auxiliary conditional gradient sequence $\seq{\Tilde{u}_k}_k$ which satisfies $\Tilde{u}_0=u_0$ as well as
\begin{equation*}
    \Tilde{u}_k= u_k+ \sigma_k(v_k-u_k) \quad \text{where} \quad v_k \in \argmin_{\norm{v}_\U \leq M_0} \left \lbrack \dual{\nabla f(u_k)}{v}+g(v) \right \rbrack
\end{equation*}
and $\sigma_k \in [0,1]$ are suitable stepsizes.
Note that
\begin{equation*}
    \Psi(u_k)= \dual{\nabla f(u_k)}{u_k-v_k}+g(u_k)-g(v_k)
\end{equation*}
by construction.
\begin{theorem} \label{thm:convexsublin}
Assume that $f$ is convex and assume that $\seq{u_k}_k$ is generated by \Cref{alg:prox}. Then there holds
\begin{equation}\label{eq:convexsublin}
    r_J(u_k) \leq \frac{r_J(u_0)}{1+qk} \qquad\forall k \in \N,
\end{equation}
where
\begin{equation*}
    q= \max\set{q_1, q_2}
    ,\quad
    q_1 = \frac{1}{2} \min\set*{1, \frac{r_J(u_0)}{4L_b M^2_0} }
    ,\quad
    q_2 = \frac{L_b r_J(u_0)}{8 (L_a + L_{\nabla f})^2 M_0^2}
    .
\end{equation*}
Moreover, $\seq{u_k}_k$ admits at least one weak*-convergent subsequence and every weak* accumulation point is a minimizer of \eqref{def:problem}.
\end{theorem}
\begin{proof}
    Recall that $J(u_k) \leq J (u_0)$, i.e., $u_k \in E(u_0)$ and thus $\norm{u_k}_\U \leq M_0 $ for all $k \in \N$.
    As a consequence, $\seq{u_k}_k$ admits at least one weak* convergent subsequence and it suffices to show the sublinear convergence rate in order to conclude optimality of weak* accumulation points.

    We start by showing \eqref{eq:convexsublin} with $q = q_1$.
    For this purpose, we take some arbitrary
    $\sigma \in [0,1]$
    and
    $v_k \in \dom g$
    and set
    $u_\sigma \coloneqq u_k+ \sigma (v_k-u_k)$.
    We get
    \begin{align}
        \nonumber
        r_J (u_{k+1})-r_J (u_k) 
        &=
        J(u_{k+1})-J(u_k) \leq  \mathcal{Q}_L (u_{k+1},u_k) \leq \mathcal{Q}_L (u_{\sigma},u_k) 
        \\
        \nonumber
        &=
        \sigma\dual{\nabla f(u_k)}{v_k-u_k}+ g(u_\sigma)-g(u_k)+ \frac{L_k}{2}\norm{u_\sigma-u_k}^2_\U
        \\
        \label{eq:nice_equation_and_inequalities}
        &\leq
        \sigma\left(\dual{\nabla f(u_k)}{v_k-u_k}+ g(v_k)-g(u_k) \right)+\frac{L_k\sigma^2}{2}\norm{v_k-u_k}^2_\U
        ,
    \end{align}
    where we used the convexity of $g$ in the last inequality.
    Now, we specialize our choice of $v_k$ to be a minimizer in the definition of $\Psi(u_k)$, i.e.,
    we choose $v_k \in \dom g$ such that
    \begin{equation*}
        \Psi(u_k)
        =
        \bracks*{ \dual{\nabla f(u_k)}{u_k-v_{k}}+g(u_k)-g(v_{k})}, \quad \norm{v_{k}}_\U \leq M_0. 
    \end{equation*}
    Inserting this choice into \eqref{eq:nice_equation_and_inequalities}
    yields
    \begin{equation*}
        r_J (u_{k+1})-r_J (u_k) 
        \le
        -\sigma \Psi(u_k)+\frac{L_k\sigma^2}{2}\norm{u_k-v_k}^2_\U
        \leq
        -\sigma r_J(u_k)+ 2L_b M^2_0\sigma^2
        .
    \end{equation*}
    By minimizing the right-hand side w.r.t.\ $\sigma \in [0,1]$,
    we get
    \begin{equation*}
        r_J (u_{k+1})-r_J (u_k)
        \leq
        -\frac{1}{2} \min\set*{1, \frac{r_J(u_k)}{4L_b M^2_0} } r_J(u_k)
        \leq
        -\frac{1}{2} \min\set*{1, \frac{r_J(u_0)}{4L_b M^2_0} } \frac{r_J(u_k)^2}{r_J(u_0)}.
    \end{equation*}
    Dividing both sides by $r_J(u_0)$ and using \cite[Lemma 3.1]{Dunnimplicit},
    we arrive at \eqref{eq:convexsublin} with $q = q_1$.

    Now, we check that $q = q_2$ is valid in \eqref{eq:convexsublin}.
    From \cref{lem:descentestimate} we know
    \begin{equation*}
        r_J(u_k) - r_J(u_{k+1})
        =
        J(u_k) - J(u_{k+1})
        \ge
        \frac{L_k}{2}\norm{u_{k+1} - u_k}_{\U}^2
        =
        \frac{1}{2 L_k}\norm{\nabla f(u_k) + \xi_k}_{\Uc}^2
    \end{equation*}
    and $\xi_k \in \partial g(x_{k+1})$.
    The first inequality from \cref{rem:convdual}
    and the convexity of $f$
    yield
    \begin{equation}
        \label{eq:in_the_second_proof}
        \norm{\nabla f(u_{k+1}) + \xi_k}_{\Uc}
        \ge
        \frac{1}{2 M_0} \Psi(u_{k+1})
        \ge
        \frac{1}{2 M_0} r_J(u_{k+1})
        .
    \end{equation}
    We also have
    \begin{equation*}
        \norm{\nabla f(u_{k+1}) + \xi_k}_{\Uc}
        \le
        L_{\nabla f} \norm{u_{k+1} - u_k}_{\U}
        +
        \norm{\nabla f(u_{k}) + \xi_k}_{\Uc}
        =
        \parens*{1 + \frac{L_{\nabla f}}{L_k}}
        \norm{\nabla f(u_{k}) + \xi_k}_{\Uc}
        ,
    \end{equation*}
    see \cref{lem:descentestimate}.
    By putting these estimates together, we get
    \begin{equation*}
        r_J(u_k) - r_J(u_{k+1})
        \ge
        \frac{L_k}{8 (L_k + L_{\nabla f})^2 M_0^2} r_J(u_{k+1}).
    \end{equation*}
    From \cite[Lemma 3.1]{Dunnimplicit},
    we get \eqref{eq:convexsublin} with $q = q_2$.
\end{proof}

\section{Improved convergence rates under structural assumptions} \label{sec:improved}
In this section, we derive linear rates of convergence provided that Problem \eqref{def:problem} satisfies additional regularity and stability properties. In that regard, we first assume that $f$ is convex and \eqref{def:problem} admits a point $\bar{u} \in \Uc$ which satisfies the \textit{QG-property}, short for quadratic growth,
\begin{align} \label{eq:quadgrowthconvergence}
    \frac{\theta}{2} \norm{u-\bar{u}}^2_\U \leq J(u)-J(\bar u) \quad \forall u \in B_R(\bar u) \text{ with } J(u) \le J(\bar u) + \eta
\end{align}
for some radius $R>0$ and parameters $\theta, \eta>0$. In this case, one readily verifies that $\bar{u}$ is the unique minimizer of Problem \eqref{def:problem} and we conclude $u_k \weakstar \bar{u}$ for the whole sequence from \cref{thm:convexsublin}. If the latter can be upgraded to strong convergence on a subsequence, the following results yield  linear convergence. 
\begin{lemma}
	\label{lem:convergence_of_the_entire_sequence_under_QG}
	Let $\seq{u_k}_k$ be generated by \Cref{alg:prox}.
	If there exists a subsequence $(u_{k_j})_j$ with $u_{k_j} \to \bar u$ in $\U$ and $\bar{u}$ satisfies the QG-property \eqref{eq:quadgrowthconvergence},
	then $u_k \to \bar u$ in $\U$.
\end{lemma}
\begin{proof}
	Let $r \in (0,R)$ be arbitrary.
	From \cref{lem:usefulstuff} we know $J(u_k) \to J(\bar u)$.
	Further, $\norm{u_{k+1} - u_k}_\U \to 0$, see \cref{cor:convofsuccessive}.
	Consequently, we can choose $K \in \N$ such that
	$ \norm{ u_K - \bar u}_\U \le r/2$
	as well as
	\begin{equation*}
		\norm{u_{k+1} - u_k}_\U \le r/2
		\quad\text{and}\quad
		J(u_k) - J(\bar u)
		\le
		\min\set*{
			\eta, \frac{\theta r^2}{8}
		}
		\qquad\forall k \ge K.
	\end{equation*}
	Let us prove $\norm{u_k - \bar u}_\U \le r/2$ for all $k \ge K$ by induction.
	It is clear that this holds for $k = K$.
	From $\norm{u_k - \bar u}_\U \le r$
	we infer
	$\norm{u_{k+1} - \bar u}_\U \le \norm{u_{k+1} - u_k}_\U + \norm{u_k - \bar u}_\U \le r \le R$.
	Consequently, we can invoke \eqref{eq:quadgrowthconvergence} with $u = u_{k + 1}$
	and this directly implies $\norm{u_{k+1} - \bar u}_\U \le r/2$.
	This finishes the induction step.
	Since $r$ was arbitrary, this shows $u_k \to \bar u$ in $\U$.
\end{proof}
\begin{theorem} \label{thm:convexlinear}
Assume that $f$ is convex and let $\seq{u_k}_k$ be generated by \Cref{alg:prox}. If there exists a subsequence $(u_{k_j})_j$ with $u_{k_j} \to \bar u$ in $\U$ and $\bar{u}$ satisfies the QG-property \eqref{eq:quadgrowthconvergence},
	then\begin{align} \label{eq:convexlin}
  r_J(u_{k+1}) \leq \mu r_J(u_ k), \quad     \norm{u_k-\bar{u}}_\U \leq C \mu^{k/2}
\end{align}
for all $k$ large enough, where
$\mu = \max\set*{\frac{1}{2},1-\frac{\theta}{4L_b}} < 1$ and $C > 0$.
\end{theorem}
\begin{proof}
We argue similarly to the proof of \cref{thm:convexsublin} but now choose
$v_k \equiv \bar u$.
Consequently, \eqref{eq:nice_equation_and_inequalities} gives
\begin{align*}
  r_J (u_{k+1})-r_J (u_k)
  &
  \leq
  \sigma \left(\dual{\nabla f(u_k)}{\bar{u}-u_k} + g(\bar{u})-g(u_k)  \right)+ \frac{L_b\sigma^2}{2}\norm{u_k-\bar u}^2_\U
  \\&
  \le
  \sigma \left(f(\bar u) - f(u_k) + g(\bar{u})-g(u_k)  \right)+ \frac{L_b\sigma^2}{2}\norm{u_k-\bar u}^2_\U
  \\&
  =
  \sigma r_J(u_k) + \frac{L_b\sigma^2}{2}\norm{u_k-\bar u}^2_\U
  ,
\end{align*}
where the second inequality comes from the convexity of $f$.
Fix an arbitrary $\eta>0$ and note that we have
$u_k \in B_R(\bar u)$ and $J(\bar{u}) \leq J(u_k)+ \eta$
for all $k \in \N$ large enough due to $u_k \rightarrow \bar{u}$, by assumption, and $r_J(u_k) \rightarrow 0$ according to \cref{thm:convexsublin}.
Consequently, we can use \eqref{eq:quadgrowthconvergence} with $u = u_k$ for $k$ large enough
and this yields
$\norm{u_k-\bar{u}}_\U^2 \leq (2/\theta) r_J(u_k)$
for all $k\in\N$ large enough. 
Combining this with the previous inequality yields
\begin{align*}
    r_J (u_{k+1}) \leq\left( 1-\sigma+ \frac{L_b \sigma^2}{\theta} \right) r_J(u_k) \quad \forall \sigma \in [0,1].
\end{align*}
Minimizing w.r.t.\ $\sigma \in [0,1] $, we arrive at the claimed convergence result for $r_J(u_k)$. The convergence result for $\norm{u_k-\bar{u}}_\U$ then follows immediately from quadratic growth.
\end{proof}
The remainder of this section is devoted to deriving similar results for nonconvex $f$.
For this purpose, we replace quadratic growth by the stronger \textit{\KL-property}, short for \KLfull,
			\begin{equation}
				\label{eq:KLconvergence}
				J(u)-J(\bar u)
				\leq
				\frac{\lambda}{2} \inf\set{\norm{\nabla f(u)+ \xi}^2_X \given \xi \in \partial g(u) \cap X}
				\quad \forall u \in B_R(\bar u) \text{ with } J(u) \le J(\bar u) + \eta
			\end{equation}      
for constants $R, \lambda, \eta >0$ as well as $X \in \{\Uc, Y\} $.
Note that \KL-property with $X = \Uc$ is stronger than the \KL-property with $X = Y$.
The relation between quadratic growth and \KL\ will be further explored in \cref{sec:KL}. 
Now, we again start by showing that strong convergence of a subsequence implies convergence of the whole sequence if the accumulation point $\bar{u}$ has the \KL-property. For the sake of brevity, we shortly summarize our assumptions.
\begin{assumption} \label{ass:KLnew}
Let $\seq{u_k}_k$ be generated by \Cref{alg:prox}. Moreover, assume that there is a subsequence~$\seq{u_{k_j}}_j$ with~$u_{k_j} \rightarrow \bar{u}$ in~$\U$
such that $\bar u$ satisfies \eqref{eq:KLconvergence} with $X \in \set{\Uc, Y}$ for constants $R,\lambda, \eta >0$.
In case $X = Y$, we additionally assume that
the sequence of subgradients $\seq{\xi_k}_k$ from \cref{lem:descentestimate} satisfies $\xi_k \in Y$ for all $k \in \N$.
\end{assumption}

In this regard, we closely follow the strategy of \cite{MehlitznonLipschitz} with appropriate changes to account for the infinite-dimensional Banach space setting of the present manuscript.
Similarly to \cref{lem:convergence_of_the_entire_sequence_under_QG}, we can show convergence of the entire sequence of iterates.
\begin{lemma} \label{thm:conwholesequence}
    Let \cref{ass:KLnew} hold. Then we have $u_k \rightarrow \bar{u}$ in $\U$ and $J(u_k) \to J(\bar u)$.
\end{lemma}
\begin{proof}
Since~$\{J(u_k)\}_k$ is decreasing and since we have~$J(u_k) \rightarrow J(\bar{u})$, see \cref{lem:usefulstuff}, we get
\begin{equation} \label{eq:etaaux}
    J(\bar{u}) < J(u_k) \leq J(\bar{u})+ \eta 
\end{equation}
for all~$k \in \N$ large enough.
Since we know $\norm{u_{k+1} - u_k}_\U \to 0$ from \cref{cor:convofsuccessive}
and since we assumed $u_{k_j} \rightarrow \bar{u}$ in $\U$,
we can choose $K \in \N$
such that \eqref{eq:etaaux} holds for all $k \geq K$ as well as
\begin{align*}
	\alpha
	\coloneqq
	\norm{u_K-\bar u}_\U
	+
	\beta \sqrt{ J(u_K)- J(\bar{u}) }
	+
	\norm{u_K - u_{K - 1}}_\U
	\leq
	R
	,
\end{align*}
where
$\beta \coloneqq \frac{2\sqrt{2\lambda} ( L_{\nabla f}+ L_b)}{L_a}$.
Next,
we prove that
\begin{align} \label{eq:statementstoprove}
	u_k \in B_R (\bar{u})
	\qquad\forall k \ge K.
\end{align}
Note that the choice of $K$ already shows $u_K \in B_R(\bar u)$.
We proceed by induction and assume that
$k \ge K$ is given such that
$u_\ell \in B_R(\bar u)$ for all $\ell = K, \ldots, k$.
Since~$\bar{u}$ is assumed to possess the \KL-property,
we get
\begin{align} \label{eq:fullestsqrt}
	\sqrt{ J(u_\ell)-J(\Bar{u})}
	\leq
	\sqrt{\frac{\lambda}{2}} \inf\set{\norm{\nabla f(u_\ell)+ \xi}_{X} \given \xi \in \partial g(u_\ell) \cap X }
	\leq
	\sqrt{\frac{\lambda}{2}} ( L_{\nabla f}+ L_b)\norm{u_{\ell}-u_{\ell-1}}_\U,   
\end{align}
for all~$\ell=K, \dots, k$,
where the second estimate is due to \cref{lem:estfordist} (see also \cref{rem:replaceUc})
and \eqref{eq:condonL}.
Using the elementary inequality
\begin{equation*}
	2 (a^2 - b^2)
	\le
	\parens*{ \frac{a - b}{c} + c a }^2
	\qquad\forall a,b \in \R, c > 0
\end{equation*}
with $a = \sqrt{J(u_\ell) - J(\bar u)}$ and $b = \sqrt{J(u_{\ell+1}) - J(\bar u)}$
yields
\begin{equation*}
	\sqrt{2} ( J(u_\ell) - J(u_{\ell+1}) )
	\le
	(c^{-1} + c) \sqrt{J(u_\ell) - J(\bar u)}
	-
	c^{-1} \sqrt{J(u_{\ell+1}) - J(\bar u)}
	.
\end{equation*}
Now, we insert \eqref{eq:fullestsqrt}
and an inequality from \cref{lem:descentestimate}
to arrive at
\begin{equation*}
	\sqrt{L_a} \norm{u_{\ell+1} - u_{\ell}}_\U
	\le
	\frac{1}{c}
	\parens*{
		\sqrt{J(u_\ell) - J(\bar u)}
		-
		\sqrt{J(u_{\ell+1}) - J(\bar u)}
	}
	+
	c \sqrt{\frac{\lambda}{2}} ( L_{\nabla f}+ 2L_b)\norm{u_{\ell}-u_{\ell-1}}_\U
\end{equation*}
The choice
$ c \coloneqq \frac{\sqrt{L_a}}{ \sqrt{2 \lambda} ( L_{\nabla f}+ 2L_b) } $
gives
\begin{equation*}
	\norm{u_{\ell+1}-u_{\ell}}_\U
	\le
	\frac{1}{\beta} \left(\sqrt{ J(u_\ell)-J(\Bar{u})}-\sqrt{ J(u_{\ell+1})-J(\Bar{u})}  \right)
	+\norm{u_{\ell}-u_{\ell-1}}_\U - \norm{u_{\ell+1}-u_{\ell}}_\U
\end{equation*}
for all $\ell = K, \ldots, k$.
Since the right-hand side telescopes,
summing over $\ell=K, \dots,k $ yields
\begin{align*}
\sum_{\ell=K}^{k}  \norm{u_{\ell+1}-u_{\ell}}_\U
&\leq
\frac{1}{\beta} \parens*{ \sqrt{ J(u_K)-J(\Bar{u})}-\sqrt{ J(u_{k+1})-J(\Bar{u})} }
+
\norm{u_{K}-u_{K-1}}_\U - \norm{u_{k + 1}-u_{k}}_\U
\\
&\le
\frac{1}{\beta} \sqrt{ J(u_K)-J(\Bar{u})}
+
\norm{u_{K}-u_{K-1}}_\U
.
\end{align*}
Hence,
\begin{equation*}
	\norm{u_{k + 1} - \bar u}_\U
	\le
	\norm{u_K - \bar u}_\U
	+
	\sum_{\ell=K}^{k}  \norm{u_{\ell+1}-u_{\ell}}_\U
	\le
	\norm{u_K - \bar u}_\U
	+
	\frac{1}{\beta} \sqrt{ J(u_K)-J(\Bar{u})}
	+
	\norm{u_{K}-u_{K-1}}_\U
	\le
	R.
\end{equation*}
This finishes the induction step and proves \eqref{eq:statementstoprove}.
Passing to the limit $k \to \infty$ in the last estimate
shows
that
$\sum_{\ell=K}^{\infty}  \norm{u_{\ell+1}-u_{\ell}}_\U < \infty$.
Consequently, the sequence~$\seq{u_k}_k$ is Cauchy in $\U$.
Together with $u_{k_j} \to \bar u$ in $\U$
we then conclude~$u_k \rightarrow \bar{u}$ in~$\U$ as claimed.
The convergence of the objective values follows from \eqref{eq:fullestsqrt}.
\end{proof}
We finish by proving linear convergence of \Cref{alg:prox}.
\begin{theorem} \label{thm:convergenceresults}
Let \cref{ass:KLnew} hold. Then we have
    \begin{align} \label{eq:nonconvexlin}
J(u_{k+1})-J(\bar{u}) \leq \frac{1}{1+\sigma} (J(u_{k})-J(\bar{u}))
\quad \text{where} \quad
\sigma \coloneqq  \frac{L_a}{\lambda ( L_{\nabla f}+ L_b)^2}
\end{align}
as well as
\begin{equation*}
       \norm{u_k-\bar{u}}_\U \leq C \mu^k
\end{equation*}
for some $C > 0$, $\mu = \sqrt{1/(1+\sigma)} < 1$ and all $k\in\N$ large enough.
\end{theorem}
\begin{proof}
    According to \cref{thm:conwholesequence}, we have~$u_k \rightarrow \bar{u}$ in~$\U$ and $J(u_k) \to J(\bar u)$.
    This implies that the \KL-property \eqref{eq:KLconvergence} holds for $k \in \N$ large enough.
    Now, we follow the second part of the proof of \cref{thm:convexsublin},
    but replace \eqref{eq:in_the_second_proof} with the \KL-property.
    This yields
    \begin{align*}
        \bracks{J(u_k) - J(\bar u)}
        -
        \bracks{J(u_{k+1}) - J(\bar u)}
        &
        \ge
        \frac{1}{2 L_k} \norm{\nabla f(u_k) - \xi_k}_X^2
        \ge
        \frac{L_k}{2 (L_k + L_{\nabla f}^2)} \norm{\nabla f(u_{k+1}) - \xi_k}_X^2
        \\&
        \ge
        \frac{L_k}{\lambda (L_k + L_{\nabla f}^2)} \bracks{J(u_{k+1}) - J(\bar u)}
        \ge
        \sigma \bracks{J(u_{k+1}) - J(\bar u)}
        .
    \end{align*}
By reordering terms, we arrive at
\eqref{eq:nonconvexlin}
for all~$k\in\N$ large enough. Thus, it remains to verify the convergence of the iterates. For this purpose, we again invoke \cref{lem:descentestimate}, to obtain
\begin{align*}
    \frac{L_a}{2} \norm{u_{k+1}-u_k}^2_\U \leq J(u_k)- J(u_{k+1}).
\end{align*}
Due to \eqref{eq:nonconvexlin}, this implies
\begin{align*}
    \norm{u_{k+1}-u_k}_\U \leq C \mu^k \quad \text{for some} \quad C >0,~\mu \coloneqq \sqrt{1/(1+\sigma)}
\end{align*}
and all~$k \in \N$ large enough.
Thus, for all~$k$ large enough and all~$K \geq k$, we have
\begin{align*}
    \norm{u_{K+1}-u_k}_\U \leq \sum^{K}_{\ell=k} \norm{u_{\ell+1}-u_\ell}_\U \leq C \mu^k \sum^{\infty}_{\ell=0} \mu^\ell=  \frac{C}{1-\mu} \mu^k.    
\end{align*}
Taking the limit for~$K \rightarrow \infty$ on the left-hand side, we arrive at the claimed statement.
\end{proof}

\section{Stability and regularity properties}
\label{sec:KL}
In this section, we discuss
the quadratic growth condition
and the \KL-property,
which were important in \cref{sec:improved}
for the derivation of convergence results.
For convenience,
we repeat their definitions
and add a third property, strong metric subregularity.

\begin{definition} \label{def:KLcondition}
	Let $\bar u \in \dom g$ be given.
	\def\lhs{J(u)-J(\bar u)}%
	\def\rhs{\frac{\lambda}{2} \inf\set{\norm{\nabla f(u)+ \xi}^2_X \given \xi \in \partial g(u) \cap X }}%
	\begin{enumerate}
		\item
			We say that~$\bar{u}$ satisfies the~\textit{\KL-property}
			(short for \KLfull)
			with $X \in \set{Y, \Uc}$
			if there exist constants $R, \lambda, \eta >0$ such that
			\begin{equation*}
				\label{eq:KL}
				\mrep[r]{
					J(u)-J(\bar u)
				}{\lhs}
				\leq
				\mrep{
					\frac{\lambda}{2} \inf\set{\norm{\nabla f(u)+ \xi}^2_X \given \xi \in \partial g(u) \cap X }
				}{\rhs}
				\quad \forall u \in B_R(\bar u) \text{ with } J(u) \le J(\bar u) + \eta
				. \tag{\KL}
			\end{equation*}
		\item
			We say that~$\bar{u}$ satisfies
			quadratic growth
			with $X \in \set{Y, \Uc}$
			if there exist constants $R, \theta, \eta > 0$ such that
			\begin{equation*}
				\label{eq:QG}
				\mrep[r]{
					\frac{\theta}{2} \norm{ u - \bar u }_{\U}^2
				}{\lhs}
				\le
				\mrep{
					J(u) - J(\bar u)
				}{\rhs}
				\quad \forall u \in B_R(\bar u) \text{ with } J(u) \le J(\bar u) + \eta
				. \tag{QG}
			\end{equation*}
		\item
			We say that~$\bar{u}$ satisfies
			strong metric subregularity
			with $X \in \set{Y, \Uc}$
			if there exist constants $R, \nu, \eta > 0$ such that
			\begin{equation*}
				\label{eq:SMS}
				\mrep[r]{
					\norm{ u - \bar u }_\U
				}{\lhs}
				\le
				\mrep{
					\frac 1 \nu \inf\set{ \norm{ \nabla f(u) + \xi }_{X} \given \xi \in \partial g(u) \cap X }
				}{\rhs}
				\quad
				\forall u \in B_R(\bar u) \text{ with } J(u) \le J(\bar u) + \eta
				.\tag{SMS}
			\end{equation*}
	\end{enumerate}
\end{definition}

Due to $Y \hookrightarrow \Uc$,
the \KL-property with $X = \Uc$
implies the \KL-property with $X = Y$.
The same comment applies to strong metric subregularity.

In the following,
we investigate the relations between the three conditions from \cref{def:KLcondition}.
Our main focus lies on providing sufficient conditions
for the \KL-property since this condition was needed in \cref{thm:convergenceresults},
see \cref{ass:KLnew}.
In particular, it would be nice if the \KL-property follows from quadratic growth,
since sufficient conditions for quadratic growth are well known.

First, we show that the \KL-property follows from quadratic growth
if the growth constant $\theta$ is large enough
compared with the non-convexity of $f$.
Of course, this is rather restrictive.
\begin{proposition}
	\label{prop:quadimpliesKL}
	We assume that $\nabla f$ satisfies a one-sided Lipschitz estimate
	\begin{equation}
		\label{eq:one_sided_Lipschitz}
		\dual{
			\nabla f(u_1) - \nabla f(u_2)
		}{u_1 - u_2 }
		\ge
		-L \norm{ u_1 - u_2 }_\U^2
		\qquad
		\forall u_1, u_2 \in \dom g
	\end{equation}
	for some $L \ge 0$.
	Further, we assume that $\bar u \in \U$ satisfies \eqref{eq:QG}
	with constants $R, \theta, \eta > 0$.
	If $\theta > L$,
	then \eqref{eq:KL}
	is satisfied with
	$X = \Uc$,
	$R = R$,
	$\lambda = 4 \theta/(\theta - L)^2$,
	and $\eta = \eta$.
\end{proposition}
\begin{proof}
	Let~$u \in B_R(\bar u)$ with $J(u) \le J(\bar u) + \eta$
	as well as~$\xi \in \partial g(u)$ be given.
	From \eqref{eq:one_sided_Lipschitz}
	we get
	\begin{equation*}
		f(\bar u)
		\ge f(u) + \dual{\nabla f(u)}{\bar u - u}
		- \frac L2\norm{\bar u - u}_{\U}^2
	\end{equation*}
	and the convexity of $g$ yields
	\(
		g(\bar u) \ge g(u) + \dual{\xi}{\bar u - u}.
	\)
	Adding these inequalities gives
	\begin{equation*}
		J(\bar u)
		\ge
		J(u) + \dual{\nabla f(u) + \xi}{\bar u - u}
		- \frac L2\norm{\bar u - u}_{\U}^2
		.
	\end{equation*}
	Now, we employ the Fenchel--Young inequality to estimate the duality product
	to obtain
	\begin{equation*}
		J(\bar u)
		\ge
		J(u)
		- \frac{1}{2\nu} \norm{\nabla f(u)+ \xi}^2_{\Uc}
		- \frac {\nu + L}2\norm{\bar u - u}_{\U}^2
		,
	\end{equation*}
	where $\nu > 0$ is arbitrary.
	Next, we insert the quadratic growth from~\eqref{eq:QG}
	to arrive at
	\begin{equation*}
		J(\bar u)
		\ge
		J(u)
		- \frac{1}{2\nu} \norm{\nabla f(u)+ \xi}^2_{\Uc}
		- \frac {\nu + L}{\theta} \parens{ J(u) - J(\bar u)}
		.
	\end{equation*}
	By rearranging, we finally obtain
	\begin{align*}
		J(u)-J(\bar{u})
		\leq
		\frac{\theta}{2 \nu \parens{\theta - \nu - L}} \norm{\nabla f(u)+ \xi}^2_{\Uc}
		.
	\end{align*}    
	Setting~$\nu=(\theta - L) / 2$, the claimed statement follows.
\end{proof}
We mention that \eqref{eq:one_sided_Lipschitz}
is satisfied with $L = 0$
if and only if $f$ is convex.
We further note that
if $\U$ is a Hilbert space,
\eqref{eq:one_sided_Lipschitz} is equivalent to the convexity of
$f + \frac{L}{2} \norm{\cdot}^2_{\U}$
and this is sometimes called $L$-semiconvexity of $f$.

The next two examples show
that, in general,
quadratic growth does not imply the \KL-property.

\begin{example}
	\label{ex:growth_and_KL}
	Let us consider $\U = \R$ and $g \equiv 0$.
	For $t > 0$, we define $f(t)$ via
	\begin{equation*}
		f(t)
		\coloneqq
		\begin{cases}
			\frac 1 {12 \cdot 4^{i}} + \frac 1 2 (t - 2^{-i})^2
			& \text{if $t \in (2^{-i}, 3 \cdot 2^{-i-1})$ for some $i \in \N$}, \\
			\frac 1 {12 \cdot 4^{i-1}} - \frac 1 2 (t - 2^{-i+1})^2
			& \text{if $t \in (3 \cdot 2^{-i-1}, 2^{-i+1})$ for some $i \in \N$}.
		\end{cases}
	\end{equation*}
	Further, we set $f(0) = 0$ and $f(t) = f(-t)$ for $t < 0$.
	One can check that $f$ is continuously differentiable and
	for $t > 0$, the gradient is given by
	\begin{equation*}
		\nabla f(t)
		=
		\begin{cases}
			t - 2^{-i}
			& \text{if $t \in (2^{-i}, 3 \cdot 2^{-i-1})$ for some $i \in \N$}, \\
			2^{-i+1} - t
			& \text{if $t \in (3 \cdot 2^{-i-1}, 2^{-i+1})$ for some $i \in \N$}.
		\end{cases}
	\end{equation*}
	Consequently, $\nabla f$ is Lipschitz with constant $1$.
	Further,
	the quadratic growth condition \eqref{eq:QG}
	is satisfied at the global minimizer $\bar t = 0$ with
	$R = \infty$,
	$\theta = 1 / 7$,
	and $\eta = \infty$.
	However, since the point $t = 2^{-i}$ is stationary for all $i \in \N$,
	i.e., $\nabla f(2^{-i}) = 0$,
	the \KL-property cannot hold for $X = \R$ and for any parameters
	$R, \lambda, \eta > 0$.
	Similarly, strong metric subregularity cannot be satisfied.

	This shows that \cref{prop:quadimpliesKL}
	cannot hold without any restriction
	on the growth constant $\theta$ and the one-sided Lipschitz constant $L$.
\end{example}

In the next example,
we have more regularity for
the nonconvex part $f$ of the objective.

\begin{example}
	\label{ex:growth_and_KL_2}
	Let us consider $\U = \R$ and set $f(t) \coloneqq -\frac34 t^2$
	and $\bar g(t) \coloneqq t^2$.
	For $t > 0$, we define the function $g$ via
	\begin{equation*}
		g(t)
		\coloneqq
		\max\set*{
			\bar g(2^{n})
			+
			\frac{ (t - 2^{n}) (\bar g(2^{n+1}) - \bar g(2^{n})) }{ 2^{n+1} - 2^{n}}
			\given
			n \in \Z
		},
	\end{equation*}
	$g(-t) \coloneqq g(t)$ and $g(0) \coloneqq 0$.
	Note that the epigraph of $g$ is the convex hull of the points
	$\set{ (\pm 2^{n}, \bar g(2^{n})) \given n \in \Z}$,
	in particular, $g \ge \bar g$.
	Consequently,
	\begin{equation*}
		(f + g)(t) \ge (f + \bar g)(t) = \frac14 t^2
		\qquad\forall t \in \R
	\end{equation*}
	and, therefore, quadratic growth is satisfied.
	Nevertheless,
	for all $n \in \Z$,
	the points $2^{n}$ are stationary for $f + g$
	since $-\nabla f(2^{n}) \in \partial g(2^{n})$.
	Thus, the \KL-property and the strong metric subregularity are violated.

	Finally, we also mention that
	$f$ is $C^2$
	and that
	the second-order condition
	\begin{equation*}
		f''(0) h^2 + g''(0,0; h) > 0
		\qquad\forall h \in \R \setminus \set{0}
	\end{equation*}
	involving the second subderivative $g''(0,0;\cdot)$ of $g$
	at $0$ w.r.t.\ $0 \in \partial g(0)$
	is satisfied,
	since this condition is equivalent to the quadratic growth condition.
	Thus, the \KL-property and the strong metric subregularity do not follow from this second-order condition.
\end{example}

However,
we mention that the function $g$ in this example
is not twice epi-differentiable,
see \cite[Definition~13.6]{RockafellarWets1998} for the definition,
which also follows from \cref{thm:SSC_implies_SMS} below.

Before we prove that the \KL-property implies quadratic growth,
we state the brief observation that both conditions together
imply strong metric subregularity.
\begin{lemma}
	\label{lem:KL_and_QG_imply_SMS}
	Suppose that \eqref{eq:KL}
	and \eqref{eq:QG}
	are satisfied at $\bar u$ with constants $R, \lambda, \eta, \theta > 0$.
	Then, \eqref{eq:SMS} holds at $\bar x$ with constants $R$ and $\nu = \sqrt{\theta/\lambda}$.
\end{lemma}
\begin{proof}
	This follows immediately from combining \eqref{eq:KL}
	and \eqref{eq:QG}.
\end{proof}

For the next results,
we need an auxiliary lemma which is a slight generalization of
\cite[Lemma~4.3]{CorellaLe2024}.
\begin{lemma}
	\label{lem:aux_min}
	Let $h_1 \colon \U \to \bar\R^+$ and $h_2 \colon \U \to \bar\R$ be given functionals
	such that $u$ is a minimizer of
	\begin{equation*}
		h_1(\cdot)^{1/2} + h_2
	\end{equation*}
	and $h_1(u) > 0$.
	Then, $u$ also minimizes
	\begin{equation*}
		h_1 + 2 h_1(u)^{1/2} h_2
		.
	\end{equation*}
\end{lemma}
\begin{proof}
  Let $v \in \dom(h_1) \cap \dom(h_2)$ be arbitrary.
	From the elementary inequality
	\begin{equation*}
		2 a^{1/2} b^{1/2} \le a + b
		\qquad\forall a,b \in \R,
	\end{equation*}
	we get with $a = h_1(u)$ and $b = h_1(v)$
	the inequality
	\begin{equation*}
		h_1(u)^{1/2}
		-
		h_1(v)^{1/2}
		\ge
		\frac{
			h_1(u) - h_1(v)
		}{
			2 h_1(u)^{1/2}
		}
		.
	\end{equation*}
	Consequently
	\begin{equation*}
		0
		\ge
		h_1(u)^{1/2} + h_2(u)
		-
		\parens{
			h_1(v)^{1/2} + h_2(v)
		}
		\ge
		\frac{
			h_1(u)
			+
			2 h_1(u)^{1/2} h_2(u)
			-
			\parens{
				h_1(v)
				+
				2 h_1(u)^{1/2} h_2(v)
			}
		}{
			2 h_1(u)^{1/2}
		}
		.
	\end{equation*}
  Since $v \in \dom(h_1) \cap \dom(h_2)$ was arbitrary,
	this shows the optimality of $u$ for the functional
	$h_1 + 2 h_1(u)^{1/2} h_2$.
\end{proof}
Note that we only discussed the case $p = 2$,
but the result also holds for more general $p > 1$
as in \cite[Lemma~4.3]{CorellaLe2024}.

Next, we demonstrate that quadratic growth follows from the \KL-property.
Consequently, the \KL-property also implies strong metric subregularity by the previous lemma.
The first proof is inspired by the proof of \cite[Theorem~2.1]{CorellaLe2024}.
\begin{lemma}
	\label{prop:KL_implies_growth}
	Suppose that \eqref{eq:KL}
	is satisfied at $\bar u$ with $X = Y$ and with constants $R, \lambda, \eta > 0$.
	We further require that
	$\bar u$ is a strict minimizer of $J$ on $B_R(\bar u)$
	and that there exist $\mu > 0$ and $\bar R \in (0,R)$ with
	\begin{equation}
		\label{eq:estimate_on_boundary}
		\mu
		\le
		J(u) - J(\bar u)
		\qquad \forall u \in \partial B_{\bar R}(\bar u).
	\end{equation}
	Then, \eqref{eq:QG}
	holds radius $\bar R$,
	$\theta = \min\set{\eta/\bar R^2, \mu/\bar R^2, 1/\lambda}$,
	and $\eta = \infty$.
\end{lemma}
\begin{proof}
	We define
  $c \coloneqq \min\set{\eta, \mu}/\bar R^2$.
	Let $\hat u \in B_{\bar R}(\bar u)$ be given.
	In case
	$J(\hat u) - J(\bar u) \ge \frac c 2 \norm{\hat u - \bar u}_{\U}^2$,
	there is nothing to show.
	Otherwise,
	we get $\hat u \ne \bar u$ and, consequently, $J(\hat u) > J(\bar u)$.
	We choose an arbitrary $\varepsilon \in (0, c^{1/2} \norm{\hat u - \bar u}_{\U} / \sqrt{J(\hat u) - J(\bar u)} -1 )$
	with $\varepsilon < 1$
	and take
	$\zeta \in Y$
	such that
	\begin{equation*}
		\dual{\zeta}{\hat u - \bar u}
		\ge
		(1-\varepsilon / 2) \norm{\zeta}_Y \norm{\hat u - \bar u}_{\U}
		\qquad\text{and}\qquad
		\norm{\zeta}_Y
		=
		(1 + \varepsilon) \frac{\sqrt{J(\hat u) - J(\bar u)}}{\norm{\hat u - \bar u}_{\U}}
		.
	\end{equation*}
	Note that the choice of $\varepsilon$ implies that
	$\norm{\zeta}_Y \le c^{1/2}$.
	Next
	we consider the functional
	\begin{equation*}
		\hat J
		\coloneqq
		\sqrt{J(\cdot) - J(\bar u)} - \dual{\zeta}{\cdot - \bar u} + \delta_{B_{\bar R}(\bar u)}
		\colon
		\U \to (-\infty,\infty]
		.
	\end{equation*}
	First, we note that
	\begin{equation*}
		\dual{\zeta}{\hat u - \bar u}
		\ge
		(1 - \varepsilon / 2) \norm{\zeta}_Y \norm{\hat u - \bar u}_{\U}
		\ge
		(1 - \varepsilon / 2) (1 + \varepsilon) \sqrt{ J(\hat u) - J(\bar u) }
		>
		\sqrt{ J(\hat u) - J(\bar u) }
	\end{equation*}
	implies
	\begin{equation*}
		\hat J(\bar u)
		=
		0
		>
		\sqrt{J(\hat u) - J(\bar u)}
		-\dual{\zeta}{\hat u - \bar u}
		=
		\hat J(\hat u)
		.
	\end{equation*}
	Hence, $\bar u$ cannot be a minimizer of $\hat J$.
	Similarly, for any $u \in \partial B_{\bar R}(\bar u)$
	we have
	\begin{equation*}
		\hat J(u)
		\ge
		\mu^{1/2}
		-
		c^{1/2} \bar R
		\ge
		0
		=
		J(\bar u)
	\end{equation*}
	and this shows that minimizers of $\hat J$
	cannot lie on the boundary $\partial B_{\bar R}(\bar u)$.

	Note that a minimizer of $\hat J$ exists, since the functional is sequentially weak* lower semicontinuous
	and $B_{\bar R}(\bar u)$ is sequentially weak* compact.
	We denote by $y \in \U$ a minimizer of the functional $\hat J$.
	The arguments above show that $y \ne \bar u$ and
	that $y$ belongs to the interior of $B_{\bar R}(\bar u)$.
	Since $\bar u$ is a strict local minimizer,
	we have $J(y) > J(\bar u)$.

	By applying \cref{lem:aux_min}
	we get that $y$ also minimizes the functional
	\begin{equation*}
		J
		- 2 \sqrt{J(y) - J(\bar u)} \zeta
		+ \delta_{B_{\bar R}(\bar u)}
		.
	\end{equation*}
	Since $y$ does not lie on the boundary of the ball,
	we get
	\begin{equation}
		\label{eq:inclusion_328}
		\xi
		\coloneqq
		-\nabla f(y)
		+
		2
		\sqrt{J(y) - J(\bar u)}
		\zeta
		\in
		\partial g(y)
		\cap Y
		.
	\end{equation}

	In order to invoke the \KL-property,
	we use the optimality of $y$ to get
	\begin{equation*}
		0
		=
		\hat J(\bar u)
		>
		\hat J(y)
		=
		\sqrt{J(y) - J(\bar u)}
		-\dual{\zeta}{y - \bar u}
	\end{equation*}
	which implies
	\begin{equation*}
		J(y) - J(\bar u)
		\le
		\dual{\zeta}{y - \bar u}^2
		\le
		\norm{\zeta}_Y^2 \bar R^2
		\le
		c \bar R^2
		\le
		\eta.
	\end{equation*}
	Thus, we can use the \KL-property \eqref{eq:KL} with $u = y$ and $\xi$ from \eqref{eq:inclusion_328}.
	This yields
	\begin{equation*}
		J(y) - J(\bar u)
		\le
		\frac\lambda 2
		\norm[\Big]{
			2 \sqrt{J(y) - J(\bar u)} \zeta
		}_{Y}^2
		=
		2 \lambda \parens{J(y) - J(\bar u)} \norm{\zeta}_{Y}^2
		.
	\end{equation*}
	Since $J(y) > J(\bar u)$,
	this term cancels out and we are left with
	\begin{equation*}
		1
		\le
		2 \lambda \norm{\zeta}_{Y}^2
		=
		2 \lambda
		(1 + \varepsilon)^2 \parens{J(\hat u) - J(\bar u)} \norm{\hat u - \bar u}_{\U}^{-2}
		.
	\end{equation*}
	With $\varepsilon \searrow 0$, we get the desired
	\begin{equation*}
		\frac1{2 \lambda} \norm{\hat u - \bar u}_{\U}^2
		\le
		J(\hat u) - J(\bar u)
		.
	\end{equation*}
	This proves the claim.
\end{proof}
In infinite dimensions,
the condition \eqref{eq:estimate_on_boundary}
does not follow from the strict optimality
of $\bar u$,
since the boundary $\partial B_{\bar R}(\bar u)$
is not
(sequentially) (weak*) compact.
In absence of this condition,
the minimizer $y$ of $\hat J$ could lie on the boundary
and the optimality condition would
additionally contain an element from
the normal cone of $B_{\bar R}(\bar u)$
at $y$.
This induces problems in the proof
and it is not clear
whether the result of \cref{prop:KL_implies_growth}
holds in absence of \eqref{eq:estimate_on_boundary}.

We give a second proof of ``\eqref{eq:KL} implies \eqref{eq:QG}'' which does not use \eqref{eq:estimate_on_boundary}
but another property.
This proof is inspired by \cite[Theorem~2]{VanNgaiThera2007}.
\begin{lemma}
	\label{prop:KL_implies_growth_take_2}
	Suppose that \eqref{eq:KL}
	is satisfied at $\bar u$ with $X \in \set{Y, \Uc}$ and with constants $R, \lambda, \eta > 0$
	and we require that
	$\bar u$ is a strict minimizer of $J$ on $B_R(\bar u)$.

	In case $X = Y$,
	we further assume that for all $u, z \in \dom g$,
	$\alpha > 0$ and $y \in Y$ such that
	$u$ is a minimizer of the functional
	\begin{equation*}
		\dual{y}{\cdot} + g + \alpha \norm{\cdot - z}_{\U}
		,
	\end{equation*}
	there exists $\xi \in \partial g(u) \cap Y$
	and $\zeta \in \partial \norm{\cdot}(u - z) \cap Y$
	such that
	\begin{equation*}
		0 = y + \xi + \alpha \zeta.
	\end{equation*}

	Then,
	\eqref{eq:QG}
	holds with radius $\bar R = R/2$,
	$\theta = \min\set{ 1/\lambda, 2 \eta / \bar R^2}$
	and $\eta$.
\end{lemma}
\begin{proof}
	We define the auxiliary function $h \colon \U \to \bar\R^+$
	via
	\begin{equation*}
		h(u)
		\coloneqq
		\sqrt{
			J(u) - J(\bar u)
		}
		+
		\delta_{B_R(\bar u)}(u)
		\qquad \forall u \in \U
		.
	\end{equation*}
	Let $u \in B_{R/2}(\bar u) \setminus \set{\bar u}$
	with $J(u) \le J(\bar u) + \eta$
	be given.
	Note that $\inf h = 0$.
	Further, we have
	\begin{equation*}
		h(u)
		\le
		\inf h
		+
		\frac{h(u)}{(1 - \varepsilon) \norm{u - \bar u}_{\U}}
		(1 - \varepsilon) \norm{u - \bar u}_{\U}
	\end{equation*}
	for an arbitrary $\varepsilon \in (0,1)$.
	Thus,
	we can invoke the celebrated variational principle by Ekeland
	to obtain $z \in \U$ such that $h(z) \le h(u)$,
	$\norm{z - u}_{\U} \le (1 - \varepsilon) \norm{u - \bar u}_{\U}$
	and $z$ minimizes
	\begin{equation*}
		h(\cdot)
		+
		\frac{h(u)}{(1 - \varepsilon) \norm{u - \bar u}_{\U}} \norm{\cdot - z}_{\U}
		.
	\end{equation*}
	Note that $z \ne \bar u$
	and, therefore, $h(z) > 0$,
	since $\bar u$ is assumed to be a strict minimizer on $B_R(\bar u)$.
	Applying \cref{lem:aux_min} to
	\begin{equation*}
		h_1 = h(\cdot)^2 = J(\cdot) - J(\bar u) + \delta_{B_R(\bar u)}
		,\qquad
		h_2 = \frac{h(u)}{(1 - \varepsilon) \norm{u - \bar u}_{\U}} \norm{\cdot - z}_{\U}
		,
	\end{equation*}
	we get that $z$ also minimizes
	\begin{equation*}
		\parens{
			J(\cdot) - J(\bar u)
			+ \delta_{B_R(\bar u)}
		}
		+
		\frac{2 h(z) h(u)}{(1 - \varepsilon) \norm{u - \bar u}_{\U}} \norm{\cdot - z}_{\U}
		.
	\end{equation*}
	We abbreviate
	$\alpha \coloneqq 2 h(z) h(u) \parens{(1 - \varepsilon) \norm{u - \bar u}_{\U}}^{-1}$.
	Due to $\norm{z - \bar u}_{\U} < R$,
	$z$ is also a local minimizer of
	\begin{equation*}
		J(\cdot)
		+
		\alpha \norm{\cdot - z}_{\U}
		.
	\end{equation*}
	By convexity of $g$, it is easy to check that $z$ also minimizes
	\begin{equation*}
		\dual{\nabla f(z)}{\cdot}
		+
		g(\cdot)
		+
		\alpha \norm{\cdot - z}_{\U}
		.
	\end{equation*}
	In case $X = \Uc$ we apply the subdifferential sum rule,
	whereas
	in case $X = Y$ we invoke the additional assumption.
	Consequently,
	\begin{equation*}
		0
		=
		\nabla f(z)
		+
		\xi
		+
		\alpha \zeta
	\end{equation*}
	for some $\xi \in \partial g(z) \cap X$ and $\zeta \in X$ with $\norm{\zeta}_X \le 1$.
	From $h(z) \le h(u)$ we get $J(z) \le J(u) \le J(\bar u) + \eta$.
	Thus, we can apply \eqref{eq:KL}
	at $z$
	and
	this yields
	\begin{equation*}
		J(z) - J(\bar u)
		\le
		\frac{\lambda}{2}
		\norm{ \alpha \zeta }_X^2
		=
		\frac{\lambda}{2}
		\parens*{
			2 \sqrt{J(z) - J(\bar u)}
			\frac{\sqrt{J(u) - J(\bar u)}}{(1 - \varepsilon) \norm{u - \bar u}_{\U}}
		}^2
		.
	\end{equation*}
	Again we use that $z \ne \bar u$ implies $J(z) - J(\bar u) > 0$,
	therefore this term cancels out
	and we are left with
	\begin{equation*}
		\frac{1 - \varepsilon}{2 \lambda} \norm{ u - \bar u }_{\U}^2
		\le
		J(u) - J(\bar u)
		.
	\end{equation*}
	Passing to the limit $\varepsilon \searrow 0$
	yields the claim.
\end{proof}

Next,
we show that the \KL-property follows from strong metric subregularity.
This is already known in the finite-dimensional case $\U = \R^n$ from
\cite[Theorem~3.1]{LiaoDingZheng2024}
and
\cite[Proposition~2]{LiuPanWuYang2024}.
\begin{lemma}
	\label{lem:SMS_implies_KL}
	Suppose that
	\eqref{eq:SMS}
	holds with $X \in \set{Y, \Uc}$ and with some constants $\nu, \eta, R > 0$.
	Further, we assume that $\nabla f$ is one-sided Lipschitz with constant $L \ge 0$,
	see \eqref{eq:one_sided_Lipschitz}.
	Then,
	\eqref{eq:KL} holds with $X$,
	radius $R$,
	\begin{equation*}
		\lambda = \frac2\nu + \frac{L}{\nu^2},
	\end{equation*}
	and
	$\eta$.
\end{lemma}
\begin{proof}
	For convenience, we give the short proof.
	Let $u$ and $\xi$ as in \eqref{eq:KL} be given.
	The Lipschitz property of $f$ combined with the convexity of $g$ yields
	\begin{equation*}
		J(u) - J(\bar u)
		\le
		\dual{\nabla f(u) + \xi}{u - \bar u}_{\U}
		+
		\frac{L}{2} \norm{u - \bar u}_{\U}^2
		.
	\end{equation*}
	Now, we bound the duality product and we insert \eqref{eq:SMS}
	to obtain
	\begin{equation*}
		J(u) - J(\bar u)
		\le
		\parens*{
			\frac1\nu + \frac{L}{2 \nu^2}
		}
		\norm{\nabla f(u) + \xi}_{Y}^2
		.
	\end{equation*}
	This shows the claim.
\end{proof}
Since $\nabla f$ is assumed to be Lipschitz continuous with constant $L_{\nabla f}$,
we can always use $L = L_{\nabla f}$ in \cref{lem:SMS_implies_KL}.
In some situations, it might be possible to use a smaller value of $L$
and this strengthens the result.
In particular, $L = 0$ is valid in the case that $f$ is convex.

In the results above,
we have seen various relations
between the three conditions from \cref{def:KLcondition}.
In particular, (under rather weak additional assumptions)
\eqref{eq:KL} and \eqref{eq:SMS}
are equivalent and imply \eqref{eq:QG}.
However, the converse implication does, in general, not hold,
see \cref{prop:quadimpliesKL,ex:growth_and_KL,ex:growth_and_KL_2}.
This is unfortunate,
since \eqref{eq:QG}
is well understood in the literature,
while \eqref{eq:KL}
is needed in the convergence analysis of \cref{sec:improved}.

In what follows,
we show that certain second-order conditions
are sufficient for \eqref{eq:SMS}
and, consequently, for \eqref{eq:KL} via \cref{lem:SMS_implies_KL}.

First we observe that the function $g$ from \cref{ex:growth_and_KL_2}
is not twice epi-differentiable.
Consequently, one might hope that this additional assumption
is enough to get \eqref{eq:SMS} and, indeed, this is the case in finite dimensions.
This result is 
a special case of \cite[Corollary~4.1]{LiMengYang2023}
and \cite[Theorem~6.3]{MohammadiSarabi2020}.
For convenience,
we give the proof.
\begin{theorem}
	\label{thm:SSC_implies_SMS}
	Let us assume that
	\begin{enumerate}
		\item
			$\U = \R^n$ is equipped with the Euclidean norm,
		\item
			$f$ is twice differentiable at $\bar u \in \U$,
		\item
			$g$ is twice epi-differentiable at $\bar u$ w.r.t.\ $-\nabla f(\bar u) \in \partial g(\bar u)$, and
		\item
			the second-order condition
			\begin{equation*}
				f''(\bar u) v^2 + g''(\bar u, -\nabla f(\bar u); v) \ge \theta \norm{v}_{\U}^2
				\qquad
				\forall v \in \U
			\end{equation*}
			holds
			for some $\theta > 0$.
	\end{enumerate}
	Then, for all $\nu < \theta$, there exists $R > 0$ such that
	\eqref{eq:SMS}
	holds with $X = \R^n$
	and constants $R$, $\nu$, and $\eta = +\infty$.
\end{theorem}
\begin{proof}
	We define
	\begin{equation*}
		\hat\nu
		\coloneqq
		\lim_{R \searrow 0}
		\inf_{\substack{u \in B_R(\bar u) \setminus \set{\bar u} \\ \xi \in \partial g(u)}}
		\frac{\norm{\nabla f(u) + \xi}_{\U^*}}{\norm{ u - \bar u }_{\U}}
	\end{equation*}
	and we have to show that $\hat\nu \ge \theta$.
	To this end, we choose sequences
	$u_n \to \bar u$, $\xi_n \in \partial g(u_n)$
	such that $t_n \coloneqq \norm{u_n - \bar u}_{\U} \searrow 0$
	and
	\begin{equation*}
		\hat\nu = \lim_{n \to \infty} \frac{\norm{\nabla f(u_n) + \xi_n }_{\U^*}}{t_n}.
	\end{equation*}
	W.l.o.g.,
	we have
	\begin{equation*}
		\frac{u_n - \bar u}{t_n} \to v
		,
		\qquad
		\frac{\nabla f(u_n) - \nabla f(\bar u)}{t_n}
		\to
		f''(\bar u) v
		\quad\text{and}\quad
		\frac{\nabla f(u_n) + \xi_n}{t_n} \to f''(\bar u) v + \zeta
	\end{equation*}
	for some $v \in \U$, $\zeta \in \U^*$.
	Note that $\norm{v}_{\U} = 1$,
	$\norm{f''(\bar u) v + \zeta}_{\U^*} = \hat\nu$
	and due to
	\begin{equation*}
		\zeta
		\leftarrow
		\frac{\xi_n - (-\nabla f(\bar u))}{t_n}
		\in
		\frac{ \partial g(u_n)  - (-\nabla f(\bar u))}{t_n}
	\end{equation*}
	we have
	\begin{equation*}
		\zeta \in D(\partial g)(\bar u, -\nabla f(\bar u); v),
	\end{equation*}
	where $D(\partial g)$ denotes the proto-derivative of $\partial g$,
	see \cite[Section~8.H]{RockafellarWets1998}.
	Now, we utilize the identity
	\begin{equation*}
		D(\partial g)(\bar u, -\nabla f(\bar u); \cdot)
		=
		\partial \parens*{\frac12 g''(\bar u, -\nabla f(\bar u); \cdot) },
	\end{equation*}
	see \cite[13.40~Theorem and 13.30~Lemma]{RockafellarWets1998}.
	Since the second subderivative is homogeneous of degree $2$,
	this implies
	\begin{align*}
		g''(\bar u, -\nabla f(\bar u); v)
		=
		\dual{\zeta}{v}_{\U}
		&=
		\dual{f''(\bar u) v + \zeta}{v}_{\U}
		-
		f''(\bar u) v^2
		\\
		&\le
		\norm{f''(\bar u) v + \zeta}_{\U^*}
		-
		f''(\bar u) v^2
		=
		\hat\nu
		-
		f''(\bar u) v^2
		.
	\end{align*}
	The second-order condition implies
	$\hat\nu \ge \theta > 0$.
\end{proof}
One might try to
transfer this result to infinite-dimensional Hilbert spaces
by using the theory of \cite{Do1992}.
However, this seems to be very delicate
since we can only extract weakly convergent subsequences
and one might end up with $\norm{v}_{\U} = 0$.
We do not know if this is possible.

We mention that
for infinite-dimensional problems of the form
\begin{equation} \label{eq:constrainedBonnans}
	\text{Minimize} \quad F(x) \quad \text{subject to}\quad G(x) \in K
\end{equation}
with twice differentiable $F$ and $G$ and convex $K$,
one can apply \cite[Proposition~4.47 and Theorem~4.51]{BonnansShapiro2000}
to prove \eqref{eq:SMS} under presence of a certain strict CQ and second-order condition.
Note that this result does not apply to our problem class.

The next result is the main result of this section.
It shows that a certain regularity property of $g$
together with a second-order condition
implies \eqref{eq:SMS}.
Due to \cref{lem:SMS_implies_KL},
this also results in \eqref{eq:KL} being satisfied.
\begin{theorem}
	\label{thm:second_order_SMS}
	Let $\bar u \in \dom g$ be given such that
	$\bar\xi \coloneqq -\nabla f(\bar u) \in \partial g(\bar u)$.
	We additionally assume the following
	for $X \in \set{Y, \Uc}$.
	\begin{enumerate}
		\item
			\label{thm:second_order_SMS:1}
			There exists an operator $f''(\bar u) \colon \U \to Y$
			such that
			\begin{equation*}
				\frac{
					\nabla f(\bar u + t_n h_n) - \nabla f(\bar u)
				}{t_n}
				\to
				f''(\bar u) h
				\qquad\text{in $Y$}
			\end{equation*}
			for all sequences $\seq{t_n}_n \subset (0,\infty)$ 
			and $\seq{h_n}_n \subset \U$
			with
			$t_n \to 0$, $h_n \weakstar h$
			and $\bar u + t_n h_n \in \dom g$.
		\item
			\label{thm:second_order_SMS:2}
			There do not exist sequences
			$\seq{t_n}_n \subset (0,\infty)$, $\seq{h_n}_n \subset \U$ and $\seq{z_n}_n \subset X$
			with
			\begin{equation*}
				t_n \to 0
				, \qquad
				h_n \weakstar 0
				, \qquad
				\norm{h_n}_{\U} = 1
				, \qquad
				z_n \to 0
				, \qquad
				\bar\xi + t_n z_n \in \partial g( \bar u + t_n h_n )
				.
			\end{equation*}
		\item
			\label{thm:second_order_SMS:4}
			For all sequences $\seq{t_n}_n \subset (0,\infty)$, $\seq{h_n}_n \subset \U$ and $\seq{z_n}_n \subset X$,
			and $z \in Y$
			with
			\begin{equation*}
				t_n \to 0
				, \qquad
				h_n \weakstar h
				, \qquad
				\norm{h_n}_{\U} = 1
				, \qquad
				z_n \to z
				, \qquad
				\bar\xi + t_n z_n \in \partial g( \bar u + t_n h_n )
			\end{equation*}
			we have
			\begin{equation*}
				\dual{z}{h}
				=
				\lim_{n \to \infty}
				\dual{ z_n }{ h_n }
				\ge
				g''(\bar u, -\nabla f(\bar u); h)
				.
			\end{equation*}
		\item
			\label{thm:second_order_SMS:5}
			The second-order condition
			\begin{equation*}
				f''(\bar u) h^2 + g''(\bar u, -\nabla f(\bar u); h)
				>
				0
				\qquad
				\forall h \in \U \setminus \set{0}
			\end{equation*}
			holds,
			where $f''(\bar u) h^2 \coloneqq \dual{f''(\bar u) h}{h}$.
	\end{enumerate}
	Then, there exist $\nu,R > 0$
	such that
	\eqref{eq:SMS} is satisfied
	with $X$ and constants
	$R$, $\nu$ and $\eta = \infty$.
\end{theorem}
\begin{proof}
	We argue by contradiction.
	If the assertion is not satisfied,
	we get sequences
	$\seq{u_n}_n \subset \U \setminus \set{\bar u}$,
	$\seq{\xi_n}_n \subset X$
	with
	$u_n \to \bar u$,
	$\xi_n \in \partial g(u_n)$
	and
	$0 = \lim_{n \to \infty} \norm{\nabla f(u_n) + \xi_n}_{X}/\norm{u_n - \bar u}_{\U}$.
	We set $t_n \coloneqq \norm{u_n - \bar u}_{\U}$.
	W.l.o.g., we assume that $h_n \coloneqq (u_n - \bar u) / t_n \weakstar h$ in $\U$
	for some $h \in \U$.

	The differentiability assumption on $\nabla f$
	in \ref{thm:second_order_SMS:1}
	ensures
	\begin{equation}
		\label{eq:z_n_and_z}
		X
		\ni
		z_n
		\coloneqq
		\frac{\xi_n - \bar \xi}{t_n}
		=
		\frac{\nabla f(u_n) + \xi_n}{t_n}
		-
		\frac{\nabla f(\bar u + t_n h_n) - \nabla f(\bar u)}{t_n}
		\to
		0
		-
		f''(\bar u) h
		=:
		z
		\in Y,
	\end{equation}
	where the convergence takes place in $X$.
	Now, $h = 0$ cannot occur,
	since this would contradict assumption \ref{thm:second_order_SMS:2}.
	Hence, $h \ne 0$.
	Using again
	\eqref{eq:z_n_and_z}
	in combination with
	assumption \ref{thm:second_order_SMS:4}
	gives
	\begin{equation*}
		-f''(\bar u) h^2
		=
		\dual{z}{h}
		\ge
		g''(\bar u, -\nabla f(\bar u); h)
	\end{equation*}
	and this contradicts
	the second-order condition \ref{thm:second_order_SMS:5}.
	This contradiction finishes the proof.
\end{proof}
We shed some light on the assumptions of \cref{thm:second_order_SMS}.
The first assumption
\ref{thm:second_order_SMS:1}
is some strong kind of second-order differentiability of $f$.
Since our requirements on $g$ are very weak, we need this strong assumption on $f$.
The second assumption
\ref{thm:second_order_SMS:2}
prevents some degenerate situations,
similarly to
the NDC in earlier works on second-order conditions
\cite{ChristofWachsmuth2017:1,wachsmuth2,BorchardWachsmuth2024}.
As in these works, one can show that it follows from a first-order growth,
see \cref{lem:prove_almost_NDC} below.
It is also clear that
\ref{thm:second_order_SMS:2}
automatically holds if $\U$ is finite dimensional.
The third assumption
\ref{thm:second_order_SMS:4}
gives a relation between the limit of difference quotients of subdifferentials of $g$
and the second subdifferential of $g$.
This is our substitute for the proto-differentiability of $\partial g$,
which is equivalent to the second epi-differentiability of $g$ for $\U = \R^n$.
The fourth assumption
\ref{thm:second_order_SMS:5}
is the usual second-order condition.
Under suitable assumptions,
this condition is equivalent to \eqref{eq:QG},
see
\cite[Theorem~4.5]{ChristofWachsmuth2017:1},
\cite[Theorem~2.10]{wachsmuth2},
or
\cite[Theorem~2.20]{BorchardWachsmuth2024}.

Note that \cref{ex:growth_and_KL}
violates
\cref{thm:second_order_SMS}~\ref{thm:second_order_SMS:1}
while
\cref{thm:second_order_SMS}~\ref{thm:second_order_SMS:4}
does not hold in
\cref{ex:growth_and_KL_2}.
It is easy to check that all the other assumptions hold in these two examples.

As announced above,
we give a sufficient condition for
assumption
\cref{thm:second_order_SMS}~\ref{thm:second_order_SMS:2}.
This sufficient condition has been previously used to prove \eqref{eq:QG},
see
\cite{WachsmuthWalter2024}.
The case $m = 0$ 
in \eqref{eq:quad_growth_with_m}
was treated in
\cite[Lemma~5.1(iii)]{ChristofWachsmuth2017:1}
and
\cite[Lemma~2.12]{wachsmuth2}.
\begin{lemma}
	\label{lem:prove_almost_NDC}
	Let $\bar u \in \dom g$ be given such that
	$\bar\xi \coloneqq -\nabla f(\bar u) \in \partial g(\bar u)$.
	Suppose that
	there exists $\theta, r > 0$
	as well as $m \in \N$, $\zeta_i \in Y$ for $i = 1,\ldots, m$,
	such that
	\begin{equation}
		\label{eq:quad_growth_with_m}
		g(u) - g(\bar u) + \dual{\nabla f(\bar u)}{u - \bar u}
		\ge
		\frac{\theta}{2} \norm{u - \bar u}_{\U}^2
		-
		\frac12 \sum_{i = 1}^m \dual{\zeta_i}{u - \bar u}^2
		\qquad
		\forall u \in B_r(\bar u)
		.
	\end{equation}
	Then,
	\cref{thm:second_order_SMS}~\ref{thm:second_order_SMS:2}
	is satisfied
	for $X = \Uc$ and for $X = Y$.
\end{lemma}
\begin{proof}
	Let us assume the existence of a sequence as in
	\cref{thm:second_order_SMS}~\ref{thm:second_order_SMS:2}.
	We define
	$u_n \coloneqq \bar u + t_n h_n$,
	$\xi_n \coloneqq \bar \xi + t_n z_n \in \partial g(u_n)$
	and
	$\bar\xi \coloneqq -\nabla f(\bar u) \in \partial g(\bar u)$.
	The convexity of $g$ together with the assumption
	yields
	\begin{align*}
		0
		&=
		\lim_{n \to \infty}
		\dual{z_n}{h_n}
		=
		\lim_{n \to \infty}
		\dual[\bigg]{
			\frac{\xi_n - \bar \xi}{t_n}
		}{\frac{u_n - \bar u}{t_n}}
		\ge
		\limsup_{n \to \infty}
		\frac{
			g(u_n) - g(\bar u) + \dual{\nabla f(\bar u)}{u_n - \bar u}
		}{t_n^2}
		\\&
		\ge
		\frac{\theta}{2}
		-
		\frac12 \liminf_{n \to \infty} \sum_{i = 1}^m \dual{\zeta_i}{h_n}^2
		=
		\frac{\theta}{2}
	\end{align*}
	and this is a contradiction.
\end{proof}

\section{An application to bang-bang control} \label{sec:bangbang}
In this section,
we are going to apply the above theory to an optimal control problem
with an optimal control which is assumed to be of bang-bang type.
As a model problem,
we consider
\begin{equation*}
	\min_{u \in \Uad} \frac12 \norm{S(u) - y_d}_{L^2(\Omega)}^2
	=
	\min_{u \in L^\infty(\Omega)} \frac12 \norm{S(u) - y_d}_{L^2(\Omega)}^2 + \delta_{\Uad}(u)
	,
\end{equation*}
where $\Omega \subset \R^d$, $d \in \set{1,2,3}$,
is open, bounded, and possesses a Lipschitz boundary,
$y_d \in L^2(\Omega)$,
$S \colon L^2(\Omega) \to H_0^1(\Omega) \cap C_0(\Omega)$
is the solution map $u \mapsto y$ of (the weak formulation of)
\begin{equation*}
	-\Delta y + a(y) = u \quad\text{in } \Omega,
	\qquad
	y = 0 \quad\text{on } \partial\Omega,
\end{equation*}
where $a \colon \R \to \R$
is monotonically increasing and $C^2$
(more general conditions are possible),
and
the feasible set is
\begin{equation*}
	\Uad
	=
	\set{
		u \in L^2(\Omega)
		\given
		u_a \le u \le u_b
	}
\end{equation*}
for some constants $u_a < u_b$.
Existence and regularity results for the control-to-state map $S$
can be found in \cite{Casas2012}.
We define
\begin{equation*}
	f(u)
	\coloneqq
	\frac12 \norm{S(u) - y_d}_{L^2(\Omega)}^2
	.
\end{equation*}
Let $\bar u \in \Uad$ be a control
and $\bar y \coloneqq S(\bar u)$ be the associated state.
The adjoint state $\bar p \coloneqq \nabla f(\bar u)$
is the (weak) solution of the adjoint equation
\begin{equation*}
	-\Delta \bar p + a'(\bar y) \bar p = \bar y - y_d \quad\text{in }\Omega,
	\qquad
	\bar p = 0 \quad\text{on } \partial\Omega.
\end{equation*}
The usual regularity results imply $\bar p \in C_0(\Omega)$.
The first-order necessary condition reads $-\bar p \in \partial\delta_{\Uad}(\bar u) $.
In line with the notation above,
we further use $\bar\xi \coloneqq -\bar p = -\nabla f(\bar u)$.
We assume that the regularity condition
\begin{equation}
	\label{eq:nondegeneracy_p_1}
	\lambda^d\parens[\big]{
		\set{ \abs{\bar\xi} \le \varepsilon}
	}
	\le
	C \varepsilon
	\qquad\forall \varepsilon > 0
\end{equation}
is satisfied,
where
$C > 0$ is a constant.
This implies that $\bar u$ is bang-bang, i.e., it attains only the values $u_a$ and $u_b$
(up to a null set).
Note that this assumption is popular for analyzing
bang-bang controls,
cf.\ \cite{Wachsmuth2008,DeckelnickHinze2012,WachsmuthWachsmuth2009,vonDaniels2016}.

As seen from the references
\cite{ChristofWachsmuth2017:1,wachsmuth2},
it is convenient to analyze the problem in the space
$\U \coloneqq \MM(\Omega)$
which is isometrically isomorphic to the dual space of
$Y \coloneqq C_0(\Omega)$.
In the sequel,
a measure from $\MM(\Omega)$ which is absolutely continuous
w.r.t.\ the Lebesgue measure will be identified with its density,
which belongs to $L^1(\Omega)$.
Under this identification, $\Uad$
will be viewed as a subset of $\MM(\Omega)$.
Next,
we are going to extend the solution operator $S$ to an operator $S \colon \MM(\Omega) \to L^2(\Omega)$.
This, however, requires some growth conditions on $a$,
in particular, $a(t) \le C (1 + \abs{t}^r)$ for some $r \in (0,3)$,
see the discussion in \cite[Section~2]{CasasKunisch2014}.
Let us argue that this is not restrictive in our situation.
Due to the control constraints,
one can check that $\norm{S(u)}_{L^\infty(\Omega)} \le C$
for all $u \in \Uad$
using the classical arguments due to Stampacchia.
Consequently, only the values of $a$ on the compact interval $[-C, C]$
are of interest and we modify $a$ outside of this interval
in order to achieve that $a$, $a'$ and $a''$ are bounded on all of $\R$.
Thus, \cite[Theorem~2.1]{CasasKunisch2014} is applicable
and this yields
that the operator $S$ can be extended to an operator $S \colon \MM(\Omega) \to L^2(\Omega)$,
since the dimension satisfies $d \le 3$.
Consequently, the objective function $f$ from above can be extended to a function
$f \colon \MM(\Omega) \to \R^+$.
Therefore, our control problem can be written as
\begin{equation*}
	\min_{u \in \MM(\Omega)} f(u) + g(u),
\end{equation*}
where
$g = \delta_{\Uad}$
and we use the above identification to consider $\Uad$
as a subset of $\MM(\Omega)$.
Note that all
$u \in \dom g = \Uad$
belong to $L^\infty(\Omega)$
and
\begin{equation*}
  \norm{u}_{\U}
  =
  \norm{u}_{\MM(\Omega)}
  =
  \norm{u}_{L^1(\Omega)}
\end{equation*}
in this case.

By definition, for $u \in \MM(\Omega)$,
the convex subdifferential $\partial g(u)$ belongs to the space
$\MM(\Omega)^*$.
As usual, we identify $C_0(\Omega)$
with a subspace of $\MM(\Omega)^* = C_0(\Omega)^{**}$
(via the canonical embedding)
and
consider
\begin{align*}
	\partial g(u) \cap C_0(\Omega)
	&=
	\set{
		\varphi \in C_0(\Omega)
		\given
		\forall v \in \MM(\Omega) :
		g(v) \ge g(u) + \dual{\varphi}{v - u}
	}
	\\&=
	\set{
		\varphi \in C_0(\Omega)
		\given
		\forall v \in \Uad :
		0 \ge \dual{\varphi}{v - u}
	}
	,
\end{align*}
where we used $g = \delta_{\Uad}$.

It is clear that
\cref{ass:setup}~\ref{ass:setup:1}, \ref{ass:setup:2} and \ref{ass:setup:4}
hold.
In order to apply our results from the previous sections,
we have to check additional assumptions.
In particular, we are going to verify
\cref{ass:setup}~\ref{ass:setup:3},
the additional assumption from \cref{prop:KL_implies_growth_take_2} in case $X = Y$,
and
\cref{thm:second_order_SMS}~\ref{thm:second_order_SMS:1}--\ref{thm:second_order_SMS:4}.
This will be addressed in the next two subsections.
Once these assumptions have been established,
we are able to prove the main theoretical result of this section,
which is stated as \cref{thm:equivalency_for_optimal_control_problems}.
To this end, we have to strengthen \eqref{eq:nondegeneracy_p_1}, see also
\cite[Assumptions~4.12]{wachsmuth2}.
\begin{assumption}
	\label{asm:regularity_assumptions}
	Let a feasible control $\bar u \in \Uad$ be given
	such that $\bar\xi \coloneqq -\nabla f(\bar u) \in C_0(\Omega)$
	satisfies the first-order condition
	$\bar\xi \in \partial g(\bar u)$.
	We assume that \eqref{eq:nondegeneracy_p_1} is satisfied.
	Further, we assume
	that $\bar\xi \in C^1(\Omega)$
	and
	that $\nabla \bar\xi(x) \ne 0$ for all $x \in \Omega$
	with $\bar\xi(x) = 0$.
\end{assumption}
Note that if $y_d \in L^p(\Omega)$ for $p > d$
and if $\Omega$ is sufficiently regular,
we get $\bar\xi = -\nabla f(\bar u) \in W^{2,p}(\Omega) \hookrightarrow C^1(\bar\Omega)$.

Further,
we emphasize
that the function $\bar\xi \in C_0(\Omega)$ is defined on the set $\Omega$.
Consequently,
\begin{equation*}
	\set{ \bar\xi = 0 }
	\coloneqq
	\set{ x \in \Omega \given \bar\xi(x) = 0}
\end{equation*}
is a subset of $\Omega$.
The same comment applies to similarly defined sets.
We also mention that \eqref{eq:nondegeneracy_p_1}
holds if $\bar\xi$ can be extended to a function in $C^1(\bar\Omega)$
with
$\nabla\bar\xi(x) \ne 0$ for all $x \in \bar\Omega$ with $\bar\xi(x) = 0$
and $\Omega$ satisfies some regularity assumptions,
see \cite[Lemma~3.2]{DeckelnickHinze2012}.
\begin{theorem}
	\label{thm:equivalency_for_optimal_control_problems}
	Let \cref{asm:regularity_assumptions}
	be satisfied.
	The following conditions are equivalent.
	\begin{enumerate}
		\item\label{item:QG}
			\eqref{eq:QG} with constants $R, \theta, \eta > 0$,
		\item\label{item:KL}
			\eqref{eq:KL} with $X = C_0(\Omega)$ and constants $R, \lambda, \eta > 0$,
			and $\bar u$ is a strict local minimizer,
		\item\label{item:SMS}
			\eqref{eq:SMS} with $X = C_0(\Omega)$ and constants $R, \nu, \eta > 0$,
			and $\bar u$ is a local minimizer,
		\item\label{item:SSC}
			the sufficient second-order condition
			\begin{equation*}
				\label{eq:SSC}
				f''(\bar u) v^2 + g''(\bar u, \bar \xi; v) > 0
				\qquad\forall v \in \MM(\Omega)
				.
				\tag{SSC}
			\end{equation*}
	\end{enumerate}
\end{theorem}
\begin{proof}
	First, we mention that
	the missing \cref{ass:setup}~\ref{ass:setup:3}
	will be shown in \cref{lem:another_assumption} below.
	Consequently, our findings from \cref{sec:KL}
	are applicable.

	``\ref{item:QG}$\Leftrightarrow$\ref{item:SSC}'':
	This follows as in
	\cite[Theorem~6.12]{ChristofWachsmuth2017:1} 
	or
	\cite[Section~5.1]{wachsmuth2}.

	``\ref{item:SSC}$\Rightarrow$\ref{item:SMS}'':
	The assumptions of \cref{thm:second_order_SMS}
	will be verified in
	\cref{lem:some_assumption,lem:condition_2_of_some_theorem,thm:second_convergence_towards_something} below.
	Consequently, \eqref{eq:SMS} holds.
	Since \ref{item:SSC} also implies \ref{item:QG},
	we also get that $\bar u$ is a local minimizer.
	This yields the claim.

	``\ref{item:SMS}$\Rightarrow$\ref{item:KL}'':
	From \cref{lem:SMS_implies_KL} we get that \eqref{eq:KL} is satisfied.
	Further, \eqref{eq:SMS} implies that $\bar u$ is an isolated stationary point.
	Since $\bar u$ is assumed to be a minimizer, it follows that $\bar u$ is a strict local minimizer,
	cf.\ \cref{thm:existence}.

	``\ref{item:KL}$\Rightarrow$\ref{item:QG}'':
	\cref{prop:KL_implies_growth_take_2} is applicable due to \cref{lem:subgradients_predual_bangbang_lemma} below.
\end{proof}

If these conditions are satisfied,
we obtain
fast convergence of the proximal gradient method in $\U = L^1(\Omega)$.
\begin{theorem}
  \label{thm:Convergence_bang_bang}
  Assume that the sequence $\seq{u_k}_k$ is generated by \cref{alg:prox}
  and that one of the weak* accumulation points $\bar u$
  (whose existence is guaranteed by \cref{lem:usefulstuff})
  satisfies \cref{asm:regularity_assumptions}.
  Then, the entire sequence $\seq{u_k}_k$ converges strongly towards $\bar u$
  and the convergence assertions from \cref{thm:convergenceresults} hold.
\end{theorem}
\begin{proof}
  Let $\seq{u_{k_j}}_j$ be a subsequence with $u_{k_j} \weakstar \bar u$ in $L^1(\Omega)$.
  \cref{asm:regularity_assumptions} implies that $\bar u$ is bang-bang.
  Consequently,
  \cite{Balder1986}
  upgrades the weak* convergence of $\seq{u_{k_j}}_j$
  to strong convergence in $L^1(\Omega)$.
  Now, \cref{thm:convergenceresults}
  is directly applicable
  (see also \cref{lem:subgradients_predual_bangbang} below)
  and yields the claim.
\end{proof}
\subsection{Verification of assumptions on \texorpdfstring{$f$}{f}}
In this subsection,
we verify that
$f \colon \MM(\Omega) \to \R^+$
satisfies \cref{ass:setup}~\ref{ass:setup:3}
and
\cref{thm:second_order_SMS}~\ref{thm:second_order_SMS:1}.

\begin{lemma}
	\label{lem:another_assumption}
	In the above setting,
	\cref{ass:setup}~\ref{ass:setup:3} is satisfied.
\end{lemma}
\begin{proof}
	From \cite[Theorem~2.1]{CasasKunisch2014}
	we directly infer that $S \colon \MM(\Omega) \to L^2(\Omega)$
	is sequentially weak*-to-strong continuous
	and this implies the corresponding property of $f$.
	Further, \cite[Proposition~3.1]{CasasKunisch2014}
	implies that $f \colon \MM(\Omega) \to \R$
	is of class $C^2$ and the derivative $f'(u)$ can be represented by
	a continuous function $\nabla f(u) \in C_0(\Omega)$,
	since the solution map of the adjoint equation provides this regularity.

	Finally,
	\cite[Proposition~3.1]{CasasKunisch2014}
	gives the expression
	\begin{equation*}
		f''(u)[h_1, h_2]
		=
		\int_\Omega
		\bracks*{ 1 - a''(y) p} \parens*{ S'(u)h_1 } \parens*{ S'(u)h_2 }
		\d x
		,
	\end{equation*}
	for all $u, h_1, h_2 \in \MM(\Omega)$,
	where $y = S(u)$ and $p = S'(u)^* (y - y_d)$.
	Using again a Stampacchia-type argument, we get a constant $C$ such that
	\begin{align*}
		\norm{S(u)}_{L^2(\Omega)} &\le C ( 1 + \norm{u}_{\MM(\Omega)})
		,
		&
		\norm{S'(u)^* y}_{L^\infty(\Omega)} &\le C ( 1 + \norm{y}_{L^2(\Omega)})
		,
		\\
		\norm{S'(u) h}_{L^2(\Omega)} &\le C \norm{h}_{\MM(\Omega)}
	\end{align*}
	for all $u, h \in \MM(\Omega)$ and $y \in L^2(\Omega)$.
	By combining these estimates and by using the boundedness of $a''$,
	this shows that the bilinear form $f''(u)$ is uniformly bounded on
	bounded subsets of $\MM(\Omega)$.
	This implies that $\nabla f \colon \Uad \to C_0(\Omega)$ is Lipschitz continuous.
\end{proof}

\begin{lemma}
	\label{lem:some_assumption}
	The assumption \cref{thm:second_order_SMS}~\ref{thm:second_order_SMS:1}
	is satisfied in the above setting.
\end{lemma}
\begin{proof}
	We already have seen that $f$ is of class $C^2$
	and
	by taking an adjoint, we get
	\begin{equation*}
		f''(\bar u) h
		=
		S'(\bar u)^*
		\braces{
			[1 - a''(\bar y) \bar p] S'(\bar u) h
		}
		\in
		C_0(\Omega)
	\end{equation*}
	for all $h \in \MM(\Omega)$,
	where $\bar y = S(\bar u)$ is the state
	and $\bar p = S'(\bar u)^* (\bar y - y_d)$ is the adjoint.
	Since $S'(\bar u) \colon \MM(\Omega) \to L^2(\Omega)$
	is sequentially weak*-to-strong continuous
	and since $S'(\bar u)^*$ is continuous from $L^2(\Omega)$
	to $C_0(\Omega)$,
	this implies that $f''(\bar u) \colon \MM(\Omega) \to C_0(\Omega)$
	is weak*-to-strong continuous.

	Now, let sequences $\seq{t_n}_n$ and $\seq{h_n}_n$ as in 
	\cref{thm:second_order_SMS}~\ref{thm:second_order_SMS:1}
	be given.
	From $h_n \weakstar h$ in $\MM(\Omega)$ and $t_n \searrow 0$,
	we get $t_n h_n \to 0$ in $\MM(\Omega)$.
	Since $f$ is of class $C^2$,
	$\nabla f \colon \MM(\Omega) \to C_0(\Omega)$
	is of class $C^1$.
	This implies
	\begin{equation*}
		\norm*{
			\nabla f(\bar u + t_n h_n) - \nabla f(\bar u) - f''(\bar u) (t_n h_n)
		}_{C_0(\Omega)}
		=
		\oo(\norm{t_n h_n}_{\MM(\Omega)})
		=
		\oo(t_n).
	\end{equation*}
	From the sequential weak*-to-strong continuity of $f''(\bar u)$,
	we get
	\begin{equation*}
		f''(\bar u) h_n
		\to
		f''(\bar u) h
		\qquad\text{in } C_0(\Omega).
	\end{equation*}
	Combining the last two equations yields the desired
	\begin{equation*}
		\frac{\nabla f(\bar u + t_n h_n) - \nabla f(\bar u)}{t_n}
		\to
		f''(\bar u) h
		\qquad\text{in } C_0(\Omega)
	\end{equation*}
	and this finishes the proof.
\end{proof}

\subsection{Verification of assumptions on \texorpdfstring{$g = \delta_{\Uad}$}{g}}

First, we verify the additional condition from \cref{prop:KL_implies_growth_take_2}.
\begin{lemma}
	\label{lem:subgradients_predual_bangbang_lemma}
	Let $z \in \Uad$, $\alpha \ge 0$ and $\varphi \in C_0(\Omega)$ be given
	and let $u \in \Uad$ be a minimizer of the functional
	\begin{equation*}
		\varphi + \alpha \norm{\cdot - z}_{\MM(\Omega)} + \delta_{\Uad}.
	\end{equation*}
	Then, there exists
	$\xi \in \partial \delta_{\Uad}(u) \cap C_0(\Omega) = \partial g(u) \cap C_0(\Omega)$
	and
	$\zeta \in \partial \norm{\cdot}_{\MM(\Omega)}(u - z) \cap C_0(\Omega)$
	such that
	$\varphi + \alpha\zeta + \xi = 0$.
	Consequently, the additional assumption from \cref{prop:KL_implies_growth_take_2}
	in case $X = Y$ is satisfied.
\end{lemma}
\begin{proof}
	First, we utilize that $\Uad$ is actually a subset of $L^2(\Omega)$.
	Consequently, we can study the optimization problem
	in the Hilbert space $L^2(\Omega)$.
	We apply the sum rule
	in the space $L^2(\Omega)$
	and this leads to
	$ 0 \in \varphi + \alpha \hat\zeta + \hat\xi $
	with
	\begin{align*}
		\hat\xi
		&\in \partial\delta_{\Uad}(u) \cap L^2(\Omega)
		\coloneqq
		\set*{
			w \in L^2(\Omega)
			\given
			\forall v \in \Uad:
			0
			\ge
			\scp{w}{v - u}_{L^2(\Omega)}
		}
		\\
		\hat\zeta
		&\in \partial\norm{\cdot}_{L^1(\Omega)}(u - z) \cap L^2(\Omega)
		,
		\\
		&\coloneqq
		\set*{
			w \in L^2(\Omega)
			\given
			\forall v \in L^2(\Omega):
			\norm{v}_{L^1(\Omega)}
			\ge
			\norm{u - z}_{L^1(\Omega)}
			+
			\scp{w}{v - (u - z)}_{L^2(\Omega)}
		}
		.
	\end{align*}
	It is well known that these inclusions are equivalent to
	$\hat\xi, \hat\zeta \in L^2(\Omega)$ and
	\begin{subequations}
		\label{eq:signs_for_subdifferentials}
		\begin{align}
			\hat\xi &\ge 0 \quad \text{on } \set{u = u_b},
			&
			\hat\zeta &=\mathrlap{1}\qquad\qquad \text{on } \set{u > z},
			\\
			\hat\xi &= 0 \quad \text{on } \set{u_a < u < u_b},
			&
			\hat\zeta &\in\mathrlap{[-1,1]}\qquad\qquad \text{on } \set{u = z},
			\\
			\hat\xi &\le 0 \quad \text{on } \set{u_a = u},
			&
			\hat\zeta &=\mathrlap{-1}\qquad\qquad \text{on } \set{u < z}.
		\end{align}
	\end{subequations}
	In case $\alpha = 0$,
	we set $\xi \coloneqq \hat\xi = -\varphi \in C_0(\Omega)$
	and we are done.
	Thus, it is sufficient to consider $\alpha > 0$ in the following.

	We define $\zeta \coloneqq \operatorname{proj}_{[-1,1]}(-\varphi / \alpha) \in C_0(\Omega)$
	and $\xi \coloneqq -\varphi - \alpha \zeta \in C_0(\Omega)$.
	In what follows, we check that
	$\xi$ and $\zeta$ also satisfy the sign conditions in \eqref{eq:signs_for_subdifferentials},
	and this will imply
	$\xi \in \partial \delta_{\Uad}(u)$
	and
	$\zeta \in \partial\norm{\cdot}_{\MM(\Omega)}(u - z)$.
	To this end, let $x \in \Omega$ be arbitrary.

	First, we consider the case
	$u(x) > z(x)$.
	Consequently, $\hat\zeta(x) = 1$.
	Further, $u(x) > z(x) \ge u_a$ and this implies $\hat\xi(x) \ge 0$.
	Hence, $-\varphi(x)/\alpha = \hat\xi(x)/\alpha + \hat\zeta \ge 1$
	and this implies
	$\zeta(x) = 1 = \hat\zeta(x)$
	and
	$\xi(x) = \hat\xi(x)$.
	Consequently, for these $x$, the system \eqref{eq:signs_for_subdifferentials} holds.
	The case
	$u(x) < z(x)$
	is similar.

	It remains to consider
	$u(x) = z(x)$.
	If, further,
	$u_a < u(x) < u_b$,
	we get $\hat\xi(x) = 0$ and, consequently,
	$\zeta(x) = \hat\zeta(x)$
	and
	$\xi(x) = \hat\xi(x)$.
	If, otherwise,
	$u_b = u(x) = z(x)$,
	we have $\hat\xi(x) \ge 0$.
	Consequently,
	$-\varphi(x) = \hat\xi(x) + \alpha \hat\zeta \ge -\alpha$.
	Hence,
	\begin{equation*}
		\xi(x)
		=
		-\varphi - \alpha \operatorname{proj}_{[-1,1]}(-\varphi / \alpha)
		\ge
		-\alpha + \alpha
		=
		0
	\end{equation*}
	and $\zeta(x) \in [-1, 1]$ holds by construction.
	Again, \eqref{eq:signs_for_subdifferentials} is satisfied.
	The final case
	$u_a = u(x) = z(x)$ is similar.

	Since \eqref{eq:signs_for_subdifferentials} holds for $\xi$,
	we immediately get that
	$\scp{\xi}{v - u}_{L^2(\Omega)} \le 0$
	holds  for all $v \in \Uad$.
	Together with $\xi \in C_0(\Omega)$,
	this implies $\xi \in \partial \delta_{\Uad}(u) \cap C_0(\Omega)$.
	Similarly, for $\zeta$ we get
	\begin{equation*}
		\norm{u - z}_{L^1(\Omega)}
		=
		\dual{\zeta}{u - z}
		\quad\text{and}\quad
		\forall v \in \MM(\Omega):
		\norm{v}_{\MM(\Omega)}
		\ge
		\dual{\zeta}{v}
		.
	\end{equation*}
	This implies
	\begin{equation*}
		\forall u \in \MM(\Omega):
		\norm{v}_{\MM(\Omega)}
		\ge
		\norm{u - z}_{\MM(\Omega)}
		+
		\dual{\zeta}{v - (u - z)}
		,
	\end{equation*}
	i.e.,
	$\zeta \in \partial \norm{\cdot}_{\MM(\Omega)}(u - z)$.
	This yields the claim.
\end{proof}

Next, we show that the subgradients $\xi_k$ from \cref{lem:descentestimate}
can be chosen from the space $Y = C_0(\Omega)$.
This property is crucial, see \cref{rem:replaceUc,ass:KLnew}.
\begin{lemma}
  \label{lem:subgradients_predual_bangbang}
  Let $u_k \in \Uad$ be given and,
  for some $L_k > 0$,
  let $u_{k+1} \in \Uad$
  be a solution of the subproblem \eqref{eq:prox},
  i.e.,
	\begin{align*}
    u_{k+1}
    &\in
    \argmin_{v \in \MM(\Omega)} \parens*{
      \dual{\nabla f(u_k)}{v} + g(u) + \frac{L_k}{2} \norm{v - u_k}_{\MM(\Omega)}^2
    }
    \\&=
    \argmin_{v \in \Uad} \parens*{
      \dual{\nabla f(u_k)}{v} + \frac{L_k}{2} \norm{v - u_k}_{L^1(\Omega)}^2
    }
    .
	\end{align*}
	Then, there exists $\xi_k \in \partial \delta_{\Uad}(u_{k+1}) \cap C_0(\Omega) = \partial g(u_{k+1}) \cap C_0(\Omega)$
  such that
  \begin{equation*}
    -\nabla f(u_k) - \xi_k
    \in
    \frac{L_k}{2}
    \partial
    \parens*{
      \norm{\cdot - u_k}_{\MM(\Omega)}^2
    }(u_{k+1})
    =
    L_k \norm{u_{k+1} - u_k}_{L^1(\Omega)} \partial\norm{\cdot}_{\MM(\Omega)}(u_{k+1} - u_k)
    .
  \end{equation*}
\end{lemma}
\begin{proof}
	Again, we utilize that $\Uad$ is actually a subset of $L^2(\Omega)$.
	We apply the sum rule
	and the chain rule \cite[Corollary~16.72]{BauschkeCombettes2011}
	in the space $L^2(\Omega)$
	and this leads to
	\begin{equation*}
		0
		\in
		\nabla f(u_k)
		+
		\hat\xi_k
		+
		\alpha_k
		\hat\zeta_k
	\end{equation*}
	with $\alpha_k = L_k \norm{u_{k+1} - u_k}_{L^1(\Omega)} \ge 0$,
	$\hat\xi_k \in \partial\delta_{\Uad}(u_{k+1}) \cap L^2(\Omega)$,
	$\hat\zeta_k \in \partial\norm{\cdot}_{L^1(\Omega)}(u_{k + 1} - u_k) \cap L^2(\Omega)$,
	where these sets are defined as in the proof of \cref{lem:subgradients_predual_bangbang_lemma}.
	Consequently,
	$u_{k+1}$ is a minimizer of
	\begin{equation*}
		\dual{\nabla f(u_k)}{\cdot} + \alpha_k \norm{\cdot - u_k}_{\MM(\Omega)} + \delta_{\Uad}
	\end{equation*}
	on $L^2(\Omega)$.
	Due to $\Uad \subset L^2(\Omega)$,
	$u_{k+1}$ also minimizes this functional over the set $\MM(\Omega)$.
	An application of \cref{lem:subgradients_predual_bangbang_lemma} yields the claim.
\end{proof}

In the proof of \cref{lem:subgradients_predual_bangbang_lemma}, the assumption $z \in \Uad$ was crucial.
Therefore, we need to assume $u_k \in \Uad$ in \cref{lem:subgradients_predual_bangbang}.
The next example shows that the assertion of \cref{lem:subgradients_predual_bangbang}
does not hold in absence of this assumption.
\begin{example}
	\label{ex:continuous_subgradients}
	Let $u_a = -1$, $u_b = 1$ and $\Omega = (0,1)$.
	Further, we consider the point $u_k = 2 (\chi_{(0,1/2)} - \chi_{(1/2,1)})$
	with $\nabla f(u_k) \equiv 0$.
	Then, for any $L_k > 0$, the solution $u_{k+1}$ of \eqref{eq:prox}
	is given by $u_{k+1} = \chi_{(0,1/2)} - \chi_{(1/2,1)}$
	which can be verified by the subgradients
	$\zeta_k = -\chi_{(0,1/2)} + \chi_{(1/2,1)}$
	and
	$\xi_k = -L_k \norm{u_{k+1} - u_k}_{L^1(\Omega)} \zeta_k$.
	However, it is not possible to choose subgradients which are continuous,
	since $\partial\norm{\cdot}_{L^1(\Omega)}(u_{k+1} - u_k)$ is a singleton.
\end{example}

Next,
we check that the functional $g = \delta_{\Uad}$
satisfies
the conditions
\ref{thm:second_order_SMS:2} and
\ref{thm:second_order_SMS:4}
of
\cref{thm:second_order_SMS}.
The first condition
follows from the non-degeneracy \eqref{eq:nondegeneracy_p_1} of the optimal adjoint.
\begin{lemma}
	\label{lem:condition_2_of_some_theorem}
	Suppose that $\bar\xi = -\nabla f(\bar u)$ satisfies assumption \eqref{eq:nondegeneracy_p_1}.
	Then
	\cref{thm:second_order_SMS}~\ref{thm:second_order_SMS:2} holds
	with $X = Y = C_0(\Omega)$.
\end{lemma}
\begin{proof}
	As in \cite[Proposition~2.7]{CasasWachsmuthWachsmuth2017:2},
	one can check that \eqref{eq:nondegeneracy_p_1} implies
	\begin{equation*}
		\dual{ -\bar\xi }{ u - \bar u }
		=
		\dual{ \nabla f(\bar u) }{ u - \bar u }
		\ge
		\kappa \norm{u - \bar u}_{L^1(\Omega)}^2
		=
		\kappa \norm{u - \bar u}_{\MM(\Omega)}^2
		\qquad\forall u \in \Uad
	\end{equation*}
	for some constant $\kappa > 0$.
	Due to $g = \delta_{\Uad}$,
	we can apply \cref{lem:prove_almost_NDC}
	and this yields the claim.
\end{proof}

The next goal is the verification of
\cref{thm:second_order_SMS}~\ref{thm:second_order_SMS:4}.
By $\HH^{d-1}$ we denote the usual $(d-1)$-dimensional Hausdorff measure
and $\HH^{d-1} | _{\set{\bar\xi = 0}}$
is its restriction to $\set{\bar\xi = 0}$,
i.e.,
\begin{equation*}
	\HH^{d-1} | _{\set{\bar\xi = 0}} (B)
	\coloneqq
	\HH^{d-1}\parens*{
		\set{\bar\xi = 0}
		\cap
		B
	}
\end{equation*}
for all Borel sets $B \subset \R^d$.

From \cite[Theorem~4.13]{wachsmuth2},
we get a characterization of the second subderivative of $g$,
see also \cite[Theorem~6.11]{ChristofWachsmuth2017:1}.
\begin{theorem}
	\label{thm:subderivative}
	Suppose that \cref{asm:regularity_assumptions} is satisfied.
	Then, $g = \delta_{\Uad}$ is strictly twice epi-differentiable at $\bar u$ for $\bar\xi$
	and
	\begin{equation*}
		g''(\bar u, \bar\xi; v)
		=
		\begin{cases}
			\displaystyle
			\int_{\set{\bar\xi = 0}}
			\frac{\abs{\nabla \bar\xi}}{u_b - u_a}
			\parens*{\frac{\d v}{\d \HH^{d-1}|_{\set{\bar\xi = 0}}}}^2
			\d\HH^{d-1}
			& \text{if } v \ll \HH^{d-1}|_{\set{\bar\xi = 0}}
			,
			\\
			\infty
			& \text{else}
		\end{cases}
	\end{equation*}
	for all $v \in \MM(\Omega)$.
\end{theorem}

We give a rough sketch of the upcoming proof of
\cref{thm:second_order_SMS}~\ref{thm:second_order_SMS:4}.
For simplicity (and only for this sketch), we assume that $u_b = -u_a = 1$,
that $z_n \equiv z$ is independent of $n$ and has compact support in $\Omega$,
and that $\set{ \bar\xi + t z }$ is a $\lambda^d$-null set for all $t > 0$.
In the notation of
\cref{thm:second_order_SMS}~\ref{thm:second_order_SMS:4},
we have
$\bar\xi \in \partial g(\bar u)$
and
$\bar\xi + t_n z \in \partial g(\bar u + t_n h_n)$.
This implies
\begin{equation*}
	\bar u = \sign( \bar\xi )
	\quad\text{and}\quad
	\bar u + t_n h_n = \sign( \bar\xi + t_n z )
	\qquad\text{$\lambda^d$-a.e.\ on $\Omega$}
	.
\end{equation*}
Consequently,
\begin{equation*}
	\dual{\psi}{h_n}
	=
	\frac{2}{t_n}
	\int_{\set{\sign(\bar\xi) \ne \sign(\bar\xi + t_n z)}} \psi \sign(z) \d\lambda^d
	\to
	2 \int_{\set{\bar\xi = 0}} \frac{\psi z}{\abs{\nabla\bar\xi}} \d\HH^{d-1}
	\quad
	\forall \psi \in C_0(\Omega),
\end{equation*}
where we used
\cite[Theorem~3.4]{wachsmuth2}.
This shows that the weak* limit of $h_n$
satisfies
\begin{equation*}
	h
	=
	\frac{2 z}{\abs{\nabla\bar\xi}} \HH^{d-1}|_{\set{\bar\xi = 0}}
	,
\end{equation*}
i.e., $2 z / \abs{\nabla\bar\xi}$ is the density of $h$
w.r.t.\ $\HH^{d-1} | _{\set{\bar\xi = 0}}$.
Together with
\cref{thm:subderivative}, we arrive at
\begin{equation*}
	\lim_{n \to \infty} \dual{z}{h_n}
	=
	\dual{z}{h}
	=
	2 \int_{\set{\bar\xi = 0}} \frac{z^2}{\abs{\nabla\bar\xi}} \d\HH^{d-1}
	=
	\int_{\set{\bar\xi = 0}}
	\frac{\abs{\nabla\bar\xi}}{2}
	\parens*{
		\frac{2 z}{\abs{\nabla\bar\xi}}
	}^2
	\d\HH^{d-1}
	=
	g''(\bar x, \bar\xi; h)
	.
\end{equation*}
Note that the ``$2$'' appearing in these formulas is actually $u_b - u_a$.
This shows that
\cref{thm:second_order_SMS}~\ref{thm:second_order_SMS:4}
holds for this particular case.
In what follows,
we check that we can arrive at a similar statement
in absence of the above simplifications.

We start by recalling (a part of) the assertion of \cite[Theorem~3.4]{wachsmuth2}.
\begin{theorem}
	\label{thm:taylor_expansion_integral_rn}
	Let $\bar\xi \in C^1(\Omega)$ be given such that
	$\nabla \bar\xi \ne 0$ on $\set{\bar\xi = 0}$.
	Further, let $z \in C_c(\Omega)$ be given.
	For all $t > 0$, we define $\Omega_t \coloneqq \set{ \sign(\bar\xi) \ne \sign(\bar\xi + t z)}$.
	Then, for all $\psi \in C(\Omega)$ we have
	\begin{subequations}
		\label{eq:diff_int_rn}
		\begin{align}
			\label{eq:diff_int_rn_a}
			\frac1t \int_{\Omega_t} \psi \d\lambda^d &\to \int_{\set{\bar\xi=0}}\frac{\psi \abs{z}}{\abs{\nabla \bar\xi}} \d\HH^{d-1}
			,
			\\
			\label{eq:diff_int_rn_b}
			\frac1t \int_{\Omega_t} \psi \sign(z) \d\lambda^d &\to \int_{\set{\bar\xi=0}}\frac{\psi z}{\abs{\nabla \bar\xi}} \d\HH^{d-1}
		\end{align}
	\end{subequations}
	as $t \searrow 0$.
\end{theorem}
Next, we show that we can use $z \in C(\Omega)$,
if the test function $\psi$ is compactly supported.
\begin{corollary}
	\label{cor:taylor_expansion_integral_rn}
	Let $\bar\xi \in C^1(\Omega)$ be given such that
	$\nabla \bar\xi \ne 0$ on $\set{\bar\xi = 0}$.
	Further, let $z \in C(\Omega)$ be given.
	For all $t > 0$, we define $\Omega_t \coloneqq \set{ \sign(\bar\xi) \ne \sign(\bar\xi + t z)}$.
	Then, for all $\psi \in C_c(\Omega)$
	the limits from \eqref{eq:diff_int_rn} hold as $t \searrow 0$.
\end{corollary}
\begin{proof}
	We choose $\varphi \in C_c(\Omega; [0,1])$ such that $\varphi = 1$ on $\set{\psi \ne 0}$.
	Then,
	\begin{equation*}
		\frac1t \int_{\Omega_t} \psi \d\lambda^d
		=
		\frac1t \int_{\set{\sign(\bar\xi) \ne \sign(\bar\xi + t (\varphi z))}} \psi \d\lambda^d
		\to
		\int_{\set{\bar\xi = 0}} \frac{\psi (\varphi z)}{\abs{\nabla\bar\xi}} \d\HH^{d-1}
		=
		\int_{\set{\bar\xi = 0}} \frac{\psi z}{\abs{\nabla\bar\xi}} \d\HH^{d-1}
		,
	\end{equation*}
	where we applied \eqref{eq:diff_int_rn_a}
	to $z$ replaced by $\varphi z \in C_c(\Omega)$.
	This shows \eqref{eq:diff_int_rn_a}
	and \eqref{eq:diff_int_rn_b} follows similarly.
\end{proof}

The next lemma is the key ingredient to generalize to the situation $z_n \ne z$.
\begin{lemma}
	\label{lem:estimate_measure_of_xi_plus_t_z}
	Let $\bar\xi \in C^1(\Omega)$ be given such that
	$\nabla \bar\xi \ne 0$ on $\set{\bar\xi = 0}$
	and \eqref{eq:nondegeneracy_p_1} hold.
	Let
	$z \in C_0(\Omega)$
	and
	$\seq{s_t}_{t > 0} \subset [0,\infty)$ with $s_t = \oo(t)$ as $t \searrow 0$ be given.
	Then,
	we have $\lambda^d(\set{\abs{\bar\xi + t z} \le s_t}) = \oo(t)$ as $t \searrow 0$.
\end{lemma}
\begin{proof}
	Let $\varepsilon \in (0,1/2)$ be arbitrary.
	In the following, we only consider small $t$ such that $s_t \le \varepsilon^2 t$.
	We define $A_t \coloneqq \set{\abs{\bar\xi + t z} \le s_t}$,
	$B \coloneqq \set{ \abs{z} \le \varepsilon }$,
	and $C_t \coloneqq \set{ \abs{\bar\xi} \le \varepsilon t }$.
	We have
	\begin{equation*}
		A_t
		\subset
		C_t
		\cup
		(A_t \cap B)
		\cup
		(A_t \setminus (B \cup C_t))
	\end{equation*}
	and we estimate the measures of the sets on the right-hand side.
	By \eqref{eq:nondegeneracy_p_1},
	we have
	\begin{equation*}
		\lambda^d( C_t ) \le C \varepsilon t
		.
	\end{equation*}
	Further,
	\begin{equation*}
		A_t \cap B
		\subset
		\set{\abs{\bar\xi} \le s_t + t \varepsilon}
	\end{equation*}
	and \eqref{eq:nondegeneracy_p_1} implies
	\begin{equation*}
		\lambda^d(A_t \cap B)
		\le
		C ( s_t + t \varepsilon)
		\le
		C (1+\varepsilon) \varepsilon t
		.
	\end{equation*}
	It remains to consider the set
	$A_t \setminus (B \cup C_t)$.
	Let $x \in A_t \setminus (B \cup C_t)$
	be arbitrary, i.e., we have
	\begin{equation*}
		\abs{\bar\xi(x)} > \varepsilon t,
		\qquad
		\abs{t z(x)} > \varepsilon t,
		\quad\text{and}\quad
		\abs{(\bar\xi + t z)(x)} \le s_t \le \varepsilon^2 t
		.
	\end{equation*}
	Hence, the signs of $\bar\xi(x)$ and $z(x)$ differ.
	W.l.o.g., we assume that $\bar\xi(x) > 0$ and $z(x) < 0$.
	Further,
	\begin{equation*}
		(\bar\xi + t (1+\varepsilon) z)(x)
		\le
		(\bar\xi + t  z)(x) + \varepsilon t z(x)
		<
		s_t - \varepsilon^2 t
		\le
		0
	\end{equation*}
	and
	\begin{equation*}
		(\bar\xi + t (1-\varepsilon) z)(x)
		\ge
		(1-\varepsilon) (\bar\xi + t  z)(x) + \varepsilon \bar\xi(x)
		>
		-(1-\varepsilon) s_t + \varepsilon^2 t
		\ge
		0
		.
	\end{equation*}
	This shows
	\begin{equation*}
		1 = \sign(\bar\xi(x)) = \sign((\bar\xi + t (1-\varepsilon) z)(x)) \ne \sign((\bar\xi + t (1+\varepsilon) z)(x))
		= -1
		.
	\end{equation*}
	Consequently,
	\begin{equation*}
		A_t \setminus (B \cup C_t)
		\subset
		\set{\sign(\bar\xi) \ne \sign(\bar\xi + t (1 + \varepsilon) z)}
		\setminus
		\set{\sign(\bar\xi) \ne \sign(\bar\xi + t (1 - \varepsilon) z)}
		.
	\end{equation*}
	Together with
	\begin{equation*}
		\set{\sign(\bar\xi) \ne \sign(\bar\xi + t (1 - \varepsilon) z)}
		\subset
		\set{\sign(\bar\xi) \ne \sign(\bar\xi + t (1 + \varepsilon) z)}
		,
	\end{equation*}
	this implies
	\begin{equation*}
		\lambda^d(A_t \setminus (B \cup C_t))
		\le
		\lambda^d(\set{\sign(\bar\xi) \ne \sign(\bar\xi + t (1 + \varepsilon) z)})
		-
		\lambda^d(\set{\sign(\bar\xi) \ne \sign(\bar\xi + t (1 - \varepsilon) z)})
		.
	\end{equation*}
	Since $z \in C_0(\Omega)$, the complement of $B$
	has a positive distance to the boundary $\partial\Omega$.
	Consequently, there exists $\psi \in C_c(\Omega; [0,1])$
	with $\psi \equiv 1$ on $\Omega \setminus B$ and,
	in particular, $\psi \equiv 1$ on $A_t \setminus (B \cup C_t)$.
	Thus,
	\begin{align*}
		\frac1t \lambda^d(A_t \setminus (B \cup C_t))
		&
		\le
		\frac1t \int_{\set{\sign(\bar\xi) \ne \sign(\bar\xi + t (1 + \varepsilon) z)}} \psi \d\lambda^d
		-
		\frac1t \int_{\set{\sign(\bar\xi) \ne \sign(\bar\xi + t (1 - \varepsilon) z)}} \psi \d\lambda^d
		\\&
		\to
		2 \varepsilon 
		\int_{\set{\bar\xi=0}}\frac{\psi \abs{z}}{\abs{\nabla \bar\xi}} \d\HH^{d-1}
		\le
		2 \varepsilon \norm{z}_{C_0(\Omega)}
		\int_{\set{\bar\xi=0}}\frac{1}{\abs{\nabla \bar\xi}} \d\HH^{d-1}
		=: \hat C \varepsilon
		\qquad\text{as } t \searrow 0
		,
	\end{align*}
	where we applied
	formula \eqref{eq:diff_int_rn_a}
	of
	\cref{cor:taylor_expansion_integral_rn} twice,
	namely
	to $z$ replaced by $(1+\varepsilon) z$ and $(1-\varepsilon) z$, respectively.
	Note that $\hat C$ is independent of $\varepsilon$ and $\hat C < \infty$ follows from \cref{asm:regularity_assumptions} due to \cite[Lemma~3.24]{wachsmuth2}.
	In particular,
	\begin{equation*}
		\lambda^d(A_t \setminus (B \cup C_t))
		\le 2 \hat C \varepsilon t
	\end{equation*}
	for $t > 0$ small enough (depending on $\varepsilon$).
	Putting everything together,
	we get
	\begin{equation*}
		\lambda^d(A_t)
		\le
		\lambda^d(C_t)
		+
		\lambda^d(A_t \cap B)
		+
		\lambda^d(A_t \setminus (B \cup C_t))
		\le
		C \varepsilon t + C (1+\varepsilon) \varepsilon t + 2 \hat C \varepsilon t
	\end{equation*}
	for $t > 0$ small enough (depending on $\varepsilon$).
	Since $\varepsilon \in (0,1/2)$ was arbitrary,
	this shows the claim.
\end{proof}
Now we are in position to prove the assumption
from
\cref{thm:second_order_SMS}~\ref{thm:second_order_SMS:4}.
\begin{theorem}
	\label{thm:second_convergence_towards_something}
	Let \cref{asm:regularity_assumptions} be satisfied.
	Then,
	for all sequences $\seq{t_n}_n \subset (0,\infty)$, $\seq{h_n}_n \subset \MM(\Omega)$ and $\seq{z_n}_n \subset C_0(\Omega)$
	with
	\begin{equation*}
		t_n \to 0
		, \qquad
		h_n \weakstar h
		, \qquad
		z_n \to z
		, \qquad
		\bar\xi + t_n z_n \in \partial g( \bar u + t_n h_n )
	\end{equation*}
	we have
	\begin{equation*}
		\dual{h}{z}
		=
		\lim_{n \to \infty}
		\dual{ h_n }{ z_n }
		=
		g''(\bar u, -\nabla f(\bar u); h)
		.
	\end{equation*}
	In particular,
	\cref{thm:second_order_SMS}~\ref{thm:second_order_SMS:4}
	holds with $X = Y = C_0(\Omega)$.
\end{theorem}
\begin{proof}
	We consider sequences
	$\seq{t_n}_n$,
	$\seq{h_n}_n$, and
	$\seq{z_n}_n$
	as in
	the claim.
	We also define
	\begin{equation*}
		\tilde h_n
		\coloneqq
		\frac{(u_b - u_a) \sign(z)}{t_n} \chi_{\set{\sign(\bar\xi) \ne \sign(\bar\xi + t_n z)}} \in \MM(\Omega)
		.
	\end{equation*}

	From \cref{cor:taylor_expansion_integral_rn} we get
	\begin{equation}
		\label{eq:awesome_limit}
		\lim_{n \to \infty} \dual{\psi}{\tilde h_n}
		=
		(u_b - u_a)
		\int_{\set{\bar\xi = 0}} \frac{\psi z}{\abs{\nabla\bar\xi}} \d\HH^{d-1}
		\qquad\forall \psi \in C_c(\Omega)
		.
	\end{equation}
	On the set
	$\set{ \abs{\bar\xi + t_n z} > t_n \norm{z_n - z}_{C_0(\Omega)}}$
	we have
	$0 \ne \sign(\bar\xi + t_n z) = \sign(\bar\xi + t_n z_n)$.
	Consequently, the definition of $\tilde h_n$ and
	$\bar\xi + t_n z_n \in \partial g(\bar u + t_n h_n )$
	imply
	that $\tilde h_n = h_n$ on this set.
	Now, \cref{lem:estimate_measure_of_xi_plus_t_z} yields
	\begin{equation*}
		\lambda^d\parens[\big]{
			\set{ \tilde h_n \ne h_n }
		}
		\le
		\lambda^d\parens[\big]{
			\set[\big]{ \abs{\bar\xi + t_n z} \le t_n \norm{z_n - z}_{C_0(\Omega)}}
		}
		=
		\oo(t_n)
	\end{equation*}
	as $n \to \infty$.
	Together with $\norm{\tilde h_n - h_n}_{L^\infty(\Omega)} \le 2 (u_b - u_a) / t_n$,
	this shows
	\begin{equation*}
		\norm{ \tilde h_n - h_n }_{\MM(\Omega)}
		=
		\norm{ \tilde h_n - h_n }_{L^1(\Omega)}
		\le
		\lambda^d(\set{\tilde h_n \ne h_n})
		\norm{ \tilde h_n - h_n }_{L^\infty(\Omega)}
		\to 0
	\end{equation*}
	as $n \to \infty$.
	Thus,
	$h_n \weakstar h$ and \eqref{eq:awesome_limit} imply
	that
	\begin{align*}
		\dual{\psi}{h}
		=
		(u_b - u_a)
		\int_{\set{\bar\xi = 0}} \frac{\psi z}{\abs{\nabla\bar\xi}} \d\HH^{d-1}
		\qquad\forall \psi \in C_c(\Omega).
	\end{align*}
	Note that the right-hand side is a continuous functional on $C_0(\Omega)$
	due to \cref{asm:regularity_assumptions} and \cite[Lemma~3.24]{wachsmuth2}.
	Since $C_c(\Omega)$ is dense in $C_0(\Omega)$,
	this implies
	$h = \parens*{(u_b - u_a) z / \abs{\nabla \bar\xi}} \HH^{d-1} |_{\set{\bar\xi = 0}}$.
	Together with
	\cref{thm:subderivative}
	we get
	\begin{equation*}
		g''(\bar u, \bar\xi; h)
		=
		\int_{\set{\bar\xi = 0}}
		\frac{\abs{\nabla\bar\xi}}{u_b - u_a}
		\parens*{
			\frac{(u_b - u_a) z}{\abs{\nabla\bar\xi}}
		}^2
		\d\HH^{d-1}
		=
		(u_b - u_a)
		\int_{\set{\bar\xi = 0}}
		\frac{z^2}{\abs{\nabla\bar\xi}}
		\d\HH^{d-1}
		=
		\dual{z}{h}
		.
	\end{equation*}
	This shows the claim.
\end{proof}
\subsection{Some negative results} \label{subsec:BonnansShapirosadness}
Before proceeding to the numerical realization of the proximal gradient method,
we briefly mention that the bang-bang problem can be recast as an abstract
contrained minimization problem of the form \eqref{eq:constrainedBonnans} by
setting $F(u) = f(u)$, $G(u)=u$, and $K=\Uad$.
Naturally, this raises the question whether
strong metric subregularity, and thus \KL, can be obtained from the
second-order conditions proposed in
\cite[Proposition~4.47 and Theorem~4.51]{BonnansShapiro2000}. In the following, we briefly sketch that this is not the case neither for $\U=L^2(\Omega)$ nor $\U=L^1(\Omega)$ or $\U=\mathcal{M}(\Omega)$. Since the constraint mapping $G$ is the identity,
the regularity assumption in \cite[Proposition~4.47]{BonnansShapiro2000}
is automatically satisfied in both cases. Moreover, the second-order condition \cite[(4.126)]{BonnansShapiro2000}
reads
\begin{equation}
	\label{eq:BS_SSC}
	f''(\bar u) h^2 \ge \beta \norm{h}_{\U}^2
	\qquad
	\forall
	h \in
	\hat C_\eta(\bar u)
	:=
	\set*{
		h \in \mathcal{T}_{\Uad}(\bar u)
		\given
		f'(\bar u)h \le \eta \norm{h}_{L^1(\Omega)}
	}
\end{equation}
for some constants $\beta, \eta > 0$ and where $\mathcal{T}_{\Uad}(\bar u)$ denotes the tangent cone of $\Uad$ at $\bar{u}$ in the respective space. If we compare this second-order condition with \eqref{eq:SSC} above,
we see that the second subderivative $g''$ is missing. Now, in the case $\U=L^2(\Omega)$, \eqref{eq:BS_SSC} cannot hold, cf.\ \cite[End of Section~2]{Casas2012}. Next, we argue that  \eqref{eq:BS_SSC} is also impossible for $\U=L^1(\Omega)$. 
From the expression for $f''$, we infer
\begin{equation*}
	\abs{ f''(\bar u) h^2 }
	\le
	C \norm{S'(\bar u) h}_{L^2(\Omega)}^2
	\le
	C
	\parens[\bigg]{
		\sup_{\norm{\psi}_{L^2(\Omega)} \le 1} \dual{S'(\bar u)^* \psi}{h}
	}^2
	.
\end{equation*}
Moreover,
\begin{equation*}
	\norm{ S'(\bar u)^* \psi }_{C^{0,\alpha}(\Omega)}
	\le
	C \norm{ S'(\bar u)^* \psi }_{H^2(\Omega)}
	\le
	C \norm{ \psi }_{L^2(\Omega)}
	\le
	C
	.
\end{equation*}
Note that this requires $H^2$-regularity for the PDE,
but for the arguments below, it suffices to have this regularity in the interior,
so this is not restrictive.
In order to simplify the next argument, we assume that $0 \in \Omega$
and that $\bar\xi(x) = x_1$ in a neighborhood of $0$.
The arguments in the general case are similar, but a little bit more technical.
Consequently, $\bar u = u_a \chi_{\set{x_1 < 0}} + u_b \chi_{\set{x_1 \ge 0}}$
locally.
For small $\varepsilon, \delta > 0$, we consider the direction
$h = \chi_{ [-\varepsilon,0) \times [-\delta, \delta]^{d-1}  } - \chi_{ [0,\varepsilon] \times [-\delta, \delta]^{d-1}  }$.
It is clear that $h \in \mathcal{T}_{\Uad}(\bar u)$
and $\abs{f'(\bar u) h} = \abs{\dual{\bar\xi}{h}} \le \varepsilon \norm{h}_{L^1(\Omega)}$.
Hence, $h \in \hat C_\eta(\bar u)$
if $\delta > 0$ is small and $\varepsilon \in (0,\eta)$.
Due to the above Hölder regularity of $\varphi \coloneqq S'(\bar u)^* \psi$,
we have
\begin{align*}
	\dual{S'(\bar u)^* \psi}{h}
	&=
	\int_{[-\delta,\delta]^{d-1}}
	\int_{-\varepsilon}^0 \varphi(x) \d\lambda^1(x_1)
	-
	\int_0^{\varepsilon} \varphi(x) \d\lambda^1(x_1)
	\d\lambda^{d-1}(x_2,\ldots,x_d)
	\\&\le
	\int_{[-\delta,\delta]^{d-1}}
	\int_0^{\varepsilon} \abs{\varphi(x) - \varphi(-x_1,x_2,\ldots,x_d)} \d\lambda^1(x_1)
	\d\lambda^{d-1}(x_2,\ldots,x_d)
	\le
	C \delta^{d-1} \varepsilon^{1 + \alpha}
	.
\end{align*}
Moreover,
$\norm{h}_{L^1(\Omega)} = 2^d \delta^{d-1} \varepsilon$.
Consequently,
\begin{equation*}
	\abs{ f''(\bar u) h^2 }
	\le
	C
	\parens[\bigg]{
		\sup_{\norm{\psi}_{L^2(\Omega)} \le 1} \dual{S'(\bar u)^* \psi}{h}
	}^2
	\le
	C \varepsilon^{2\alpha} \norm{h}_{L^1(\Omega)}^2
	.
\end{equation*}
Since $\varepsilon \in (0,\eta)$ was arbitrary,
this shows that the second-order condition \eqref{eq:BS_SSC} cannot hold
and, subsequently,
\cite[Theorem~4.51]{BonnansShapiro2000}
is not applicable
in the space $L^1(\Omega)$.
Since $L^1(\Omega) \subset \MM(\Omega)$,
we can repeat the above arguments and this inhibits also the application in $\MM(\Omega)$.
\subsection{Discretization and numerical realization}

In this subsection,
we address the discretization of the control problem
and
how the proximal gradient method \cref{alg:prox}
can be implemented efficiently.

We use a standard discretization by the finite element method.
The domain $\Omega$ is replaced by a triangular (or tetrahedral) mesh.
The state $y$ and the adjoint $p$
are discretized by continuous and piecewise linear functions.
The control $u$ is discretized by piecewise constant functions,
which are identified with a vector from $\R^N$.
In each step of the proximal gradient method, we have to compute a solution
$u_{k+1} \in \R^N$ of
\begin{equation*}
	\argmin_{v \in [u_a, u_b]^N}
	\gamma_{k}^\top v + \frac{L_k}{2} \parens[\bigg]{\sum_{j = 1}^N a_j \abs{v_j - u_{k,j}}}^2
	,
\end{equation*}
where $a_j > 0$ is the area of triangle $j$
and the vector $\gamma_k \in \R^N$ represents the discretization of $\nabla f(u_k)$.
It is straightforward to check that $v$ is a solution of this problem
if and only if $w = a \odot v$ (entrywise product)
solves
\begin{equation*}
	\argmin_{w \in [\hat u_a, \hat u_b]}
	\hat\gamma_{k}^\top w + \frac{1}{2} \norm{ w - \hat u_k }_1^2
	,
\end{equation*}
where $\hat u_a = u_a a$, $\hat u_b = u_b a$, $\hat\gamma_{k,j} = \gamma_{k,j} / (a_j L_k)$, and $\hat u_k = a \odot u_k$.
From \cref{prop:wellposedalg},
we know that the value $\alpha_k = \norm{ w_k - \hat u_k }_1$ is the same for all solutions $w$.
Further, the optimality of a point $w_k$ can be characterized by
\begin{equation*}
	0
	\in
	\hat\gamma_k + \alpha_k \partial \norm{\cdot}_1(w_k - \hat u_k) + \partial\delta_{[\hat u_a, \hat u_b]}(w_k)
	.
\end{equation*}
If we would already know the value of $\alpha_k$, we could reconstruct solutions $w_k$ by
selecting values
$ w_{k,j} \in \bar W_{k,j}(\alpha_k)$,
such that $\norm{w_k - \hat u_k}_1 = \alpha_k$,
where
\begin{align*}
	\bar W_{k,j}(\alpha) &\coloneqq \mrep{\set{\hat u_{b,j}}          }{\mspace{99mu}} \text{if } \hat\gamma_{k,j} < -\alpha, &
	\bar W_{k,j}(\alpha) &\coloneqq \mrep{[\hat u_{k,j}, \hat u_{b,j}]}{\mspace{99mu}} \text{if } \hat\gamma_{k,j} = -\alpha < 0, \\
	\bar W_{k,j}(\alpha) &\coloneqq \mrep{\set{\hat u_{a,j}}          }{\mspace{99mu}} \text{if } \hat\gamma_{k,j} > \alpha, &
	\bar W_{k,j}(\alpha) &\coloneqq \mrep{[\hat u_{a,j}, \hat u_{k,j}]}{\mspace{99mu}} \text{if } \hat\gamma_{k,j} = \alpha > 0, \\
	\bar W_{k,j}(\alpha) &\coloneqq \mrep{\set{\hat u_{k,j}}          }{\mspace{99mu}} \text{if } \abs{\hat\gamma_{k,j}} < \alpha, &
	\bar W_{k,j}(\alpha) &\coloneqq \mrep{[\hat u_{a,j}, \hat u_{b,j}]}{\mspace{99mu}} \text{if } \hat\gamma_{k,j} = \alpha = 0.
\end{align*}
Note that the set-valued mapping
\begin{equation*}
	\alpha
	\mapsto
	N_k(\alpha)
	\coloneqq
	\set{ \norm{w - \hat u_k}_1 \given w_j \in \bar W_{k,j}(\alpha) }
\end{equation*}
is decreasing
and our optimization problem amounts to find the unique
$\alpha_k$ such that $\alpha_k \in N_k(\alpha_k)$.
For $\alpha > 0$, we have the representation
\begin{equation*}
	N_k(\alpha)
	=
	\sum_{j : \hat\gamma_{k,j} < -\alpha} \set{\hat u_{b,j} - \hat u_{k,j} }
	+
	\sum_{j : \hat\gamma_{k,j} > \alpha} \set{\hat u_{k,j} - \hat u_{a,j} }
	+
	\sum_{j : \hat\gamma_{k,j} = -\alpha} \bracks{0, \hat u_{b,j} - \hat u_{k,j} }
	+
	\sum_{j : \hat\gamma_{k,j} = \alpha} \bracks{0, \hat u_{k,j} - \hat u_{a,j} }
	.
\end{equation*}
Consequently, the optimal value $\alpha_k$ can be found by transferring the approach of quickselect,
see \cite{Hoare1961}.
That is, given a pivot element $\alpha$,
we can compute $N_k(\alpha)$ efficiently by partitioning the set of indices $J = \set{1,\ldots,N}$
according to whether $\abs{\hat\gamma_{k,j}} \le \alpha$ or $\abs{\hat\gamma_{k,j}} > \alpha$.
This immediately tells us whether $\alpha < \alpha_k$, $\alpha = \alpha_k$ or $\alpha > \alpha_k$
and we remains to iterate only over one part of $J$.
If we select the pivot element randomly from $\abs{\hat\gamma_{k,j}}$,
we arrive at an average complexity of $O(N)$.
Alternatively, the pivot element could be selected by a median-of-medians strategy
to obtain even a worst-case complexity of $O(N)$,
cf.\ \cite{BlumFloydPrattRivestTarjan1973}.

This discussion also shows that the next iterate $u_{k+1}$
can be selected such that it almost is bang-bang.
Indeed, if only $m < N$ components of $u_k$ are not at the bounds $\set{u_a, u_b}$,
then $u_{k+1}$ can be chosen such that at most $m + 1$ components are not at the bounds $\set{u_a, u_b}$.

\subsection{Numerical examples} \label{subsec:numerical}
\begin{figure}[ht]
	\centering
	\includegraphics{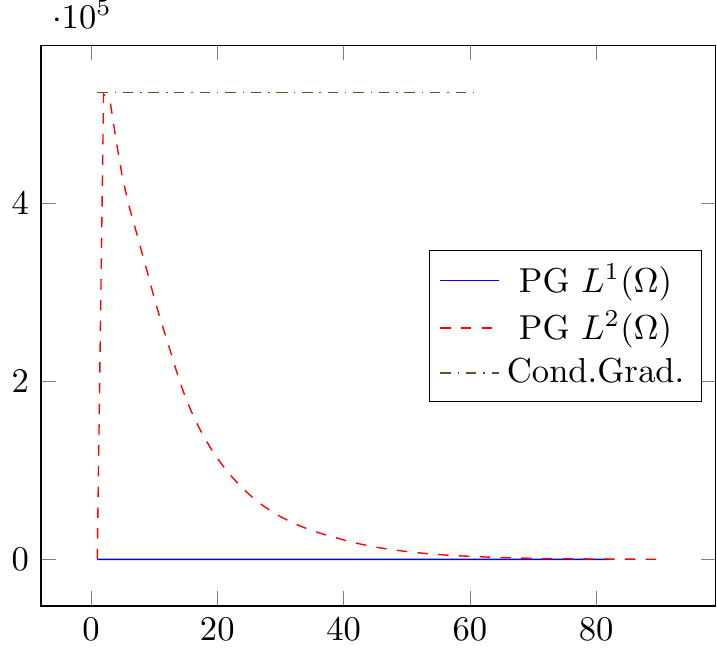}
	\includegraphics{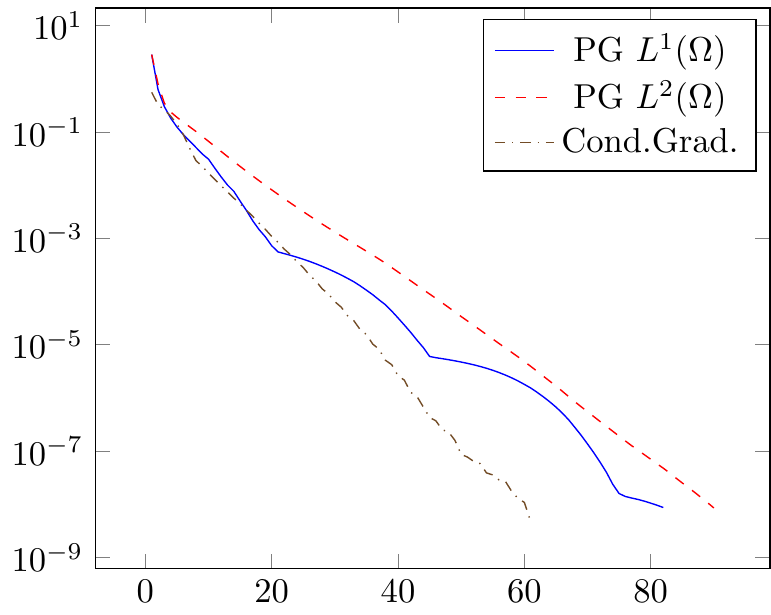}
	\caption{Number of coefficients in the control which are not bang-bang (left) and progress of the objective value (right),
	for the solution of the linear problem ($\alpha = 0$).}
	\label{fig:convergence}
\end{figure}

\begin{figure}[ht]
	\centering
	\includegraphics{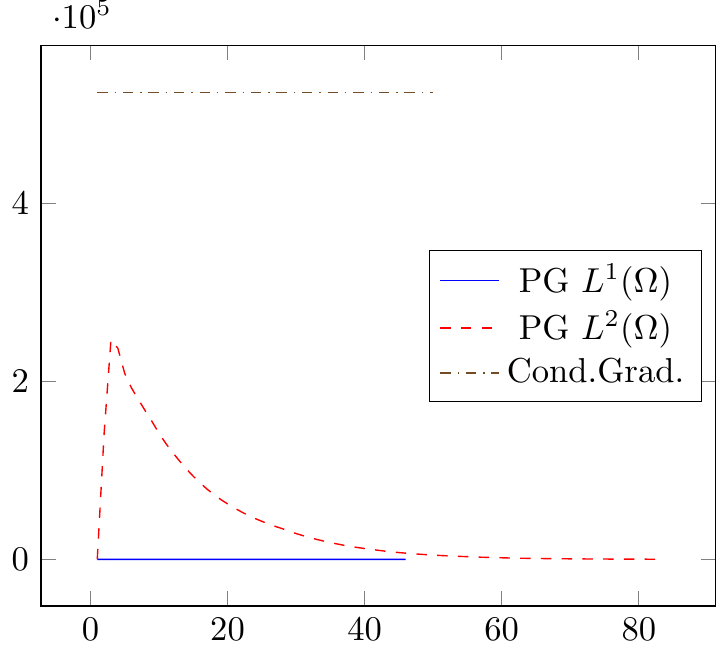}
	\includegraphics{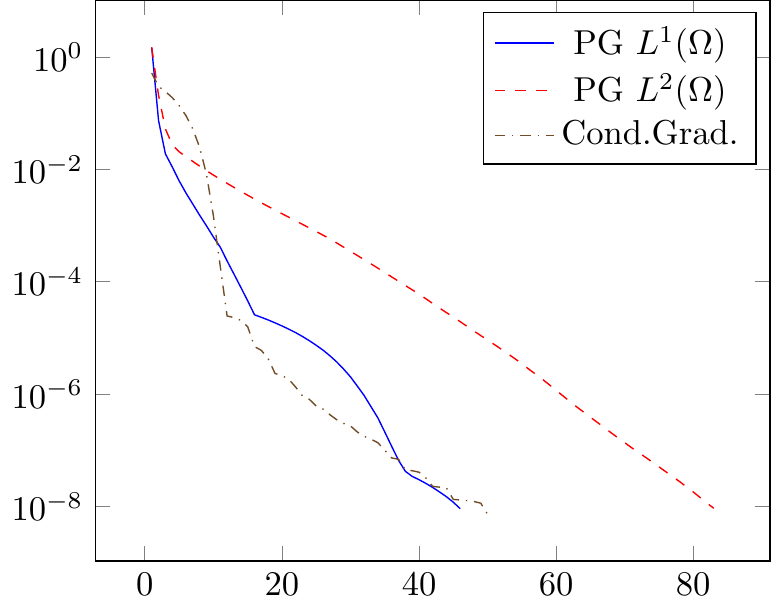}
	\caption{Number of coefficients in the control which are not bang-bang (left) and progress of the objective value (right),
	for the solution of the nonlinear problem ($\alpha = 10$).}
	\label{fig:convergence_nl}
\end{figure}

Finally, we address numerical results for the particular setting of
\begin{equation*}
	\Omega = (0,1)^2
	,\quad
	y_d(x) = 2 \sin(\pi (x_1+1)) \cos(2 \pi x_2) + 0.8 (x_2+x_1^2) - 0.5
	,\quad
	u_b = -u_a = 5
	,\quad
	a(y) = \alpha y^3
	,
\end{equation*}
where $\alpha\geq 0$. In the following, we consider a linear example, i.e., $\alpha=0$, as well as a nonlinear one choosing $\alpha=10$. In both cases, three different algorithms are applied: The proposed proximal gradient method from \cref{sec:problen} with $\U = \mathcal{M}(\Omega)$ in which we realize the subproblems as described above,
the ``usual'' proximal gradient method in $\U=L^2(\Omega)$, leading to
\begin{equation*}
    u_{k+1}= \max\left\{ \min\left\{u_k- \frac{1}{L_k} \nabla f(u_k),u_b\right\}, u_a \right\}
\end{equation*}
as well as the conditional gradient or Frank-Wolfe method
\begin{equation*}
u_{k+1}=u_k+ \sigma_k (v_k -u_k) \quad \text{where}\quad v_k(x) \in
    \begin{cases}
      \{u_a\}   & \text{if } p_k(x) >0, \\
      \{u_b\}   & \text{if } p_k(x)\le0
    \end{cases} \quad \text{for a.e. }x \in \Omega,
\end{equation*}
$p_k=\nabla f(u_k)$ and $\sigma_k \in [0,1]$.
For both proximal gradient methods, we use a backtracking line search as described in \Cref{rem:stepsizecho} 
to guarantee \eqref{eq:condonL},
while the line search from
\cite{pedregosa2020} with $\norm{\cdot}=\norm{\cdot}_{L^1(\Omega)}$ is adapted for Frank-Wolfe.
The problem is discretized by the finite element method as described above.
The proximal gradient algorithms are initialized with $u_0 = u_b$
while the Frank-Wolfe method is initialized with $u_0 = 0$.
All algorithms stopped as soon as the dual gap functional $\Psi(u_k)$ drops below $10^{-8}$.
All calculations where carried out in Matlab(R) R2024b
on a Intel(R) Core(TM) i7-8700 CPU.
The implementation is available on GitHub, see \url{https://github.com/gerw/prox_grad_banach_space}.
For both problems, we start by running all algorithms on a hierarchy of uniformly refined triangulations, see \cref{tab:linear,tab:nonlinear}.
In these tables, we show the number of iterations and the required time (in seconds) for all algorithms.
Moreover, on the mesh with $263\,169$ nodes, we provide additional convergence plots showing, on the one hand, the progression of the objective residual
$J(u_k) - J(u^*)$ as well as the, on the other hand, the number of coefficients of the iterate,
which are not bang-bang, i.e., not equal to $u_a$ or $u_b$, see \cref{fig:convergence,fig:convergence_nl}.
In this regard, $u^*$ is computed by running the proximal gradient method in $L^1(\Omega)$
with a tighter termination criterion.
\begin{table}
	\centering
	\begin{tabular}{rrrrrrrr}
		\toprule
		& & \multicolumn{2}{c}{PG $L^1(\Omega)$} & \multicolumn{2}{c}{PG $L^2(\Omega)$} & \multicolumn{2}{c}{Cond.grad.}\\
		\cmidrule(r){3-4} \cmidrule(r){5-6} \cmidrule(r){7-8}
		\#nodes &\#triang  &  iter  &  time  & iter  &  time  & iter  &  time \\
		   1089 &    2048  &   149  &  0.10  &   84  &  0.05  &  363  &  0.25 \\
		   4225 &    8192  &   144  &  0.19  &   90  &  0.11  &  275  &  0.34 \\
		  16641 &   32768  &   102  &  0.58  &   90  &  0.43  &   75  &  0.38 \\
		  66049 &  131072  &    86  &  2.36  &   89  &  2.08  &   68  &  1.86 \\
		 263169 &  524288  &    82  & 11.87  &   90  & 11.34  &   61  &  8.03 \\
		1050625 & 2097152  &    85  & 68.44  &   91  & 68.48  &   58  & 46.86 \\
		\bottomrule
	\end{tabular}
	\caption{Numerical results in the linear case}
	\label{tab:linear}
\end{table}
\begin{table}
	\centering
	\begin{tabular}{rrrrrrrr}
		\toprule
		& & \multicolumn{2}{c}{PG $L^1(\Omega)$} & \multicolumn{2}{c}{PG $L^2(\Omega)$} & \multicolumn{2}{c}{Cond.grad.}\\
		\cmidrule(r){3-4} \cmidrule(r){5-6} \cmidrule(r){7-8}
		\#nodes &\#triang  &  iter &   time  & iter  &  time  & iter  &  time  \\
		   1089 &    2048  &   131 &   0.51  &    88 &   0.35 &  303  &   1.16 \\
		   4225 &    8192  &    48 &   0.87  &    77 &   1.27 &   75  &   1.05 \\
		  16641 &   32768  &    81 &   5.84  &    84 &   6.62 &   84  &   6.53 \\
		  66049 &  131072  &    57 &  20.26  &    82 &  26.14 &   66  &  20.64 \\
		 263169 &  524288  &    46 & 101.47  &    83 & 159.16 &   50  &  99.07 \\ 
		\bottomrule
	\end{tabular}
	\caption{Numerical results in the nonlinear case}
	\label{tab:nonlinear}
\end{table}

In the linear case, and on the finest grid, we see that the performance of the three methods is comparable while $L^1$-proximal gradient as well as Frank-Wolfe outperform the $L^2$-proximal gradient method in the nonlinear case. Note that the latter exhibits mesh-independence, i.e., the number of iterations required to reach the termination criterion is almost independent of the meshwidth. In contrast, for both Frank-Wolfe and $L^1$-proximal gradient, we observe faster convergence the finer the grid is. For Frank-Wolfe, this behavior was already discussed in \cite{DunnSachs1983}. Since linear convergence of the former for linear problems can be shown under regularity assumptions akin to \eqref{eq:nondegeneracy_p_1}, see \cite{Kunisch2024}, this similarity is expected. Next, we compare the structure of the iterates in the left part of \cref{fig:convergence,fig:convergence_nl}. On the one hand, as predicted, the iterates of the proximal gradient method in $L^1$
are almost bang-bang, reflecting the expected structure of the stationary point we want to approximate. On the other hand, both, $L^2$-proximal gradient and Frank-Wolfe exhibit a significant number of non-bang-bang coefficients in intermediate iterations. In particular, while their number decreases for $L^2$-proximal gradient as $k$ grows, the same is not true for Frank-Wolfe as the stepsize typically satisfies $\sigma_k \in (0,1)$, i.e., full steps are rarely performed and $u_{k+1}$ is obtained as a true convex combination between $u_k$ and the bang-bang function $v_k$. In summary, the proposed $L^1$-proximal gradient scheme combines the convergence speed of Frank-Wolfe methods with favourably structured iterates highlighting its practical relevance.   

Finally, we repeat the experiments of \cref{tab:linear,tab:nonlinear} on a problem with Neumann boundary conditions,
i.e., the PDE is replaced by
\begin{equation*}
	-\Delta y + 10y + a(y) = u \quad\text{in } \Omega,
	\qquad
	\frac{\partial y}{\partial \nu} = 0 \quad\text{on } \partial\Omega
\end{equation*}
and for the control bounds, we use $u_b = -u_a = 10$. We emphasize that the theoretical parts of \cref{sec:bangbang}, in particular the derivation of the second subderivative, does not carry over immediately to Neumann boundary conditions. This is an open problem and beyond the scope of the current work. However, the practical realization of the $L^1$-proximal gradient works analogously.  The results are shown in \cref{tab:linear_nbc,tab:nonlinear_nbc}. For this setting, we make the same observations regarding the convergence behavior w.r.t.\ refinement of the grid, while $L^1$-proximal gradient outperforms both of its competitors on the finest discretization.
\begin{table}
	\centering
	\begin{tabular}{rrrrrrrr}
		\toprule
		& & \multicolumn{2}{c}{PG $L^1(\Omega)$} & \multicolumn{2}{c}{PG $L^2(\Omega)$} & \multicolumn{2}{c}{Cond.grad.}\\
		\cmidrule(r){3-4} \cmidrule(r){5-6} \cmidrule(r){7-8}
		\#nodes &\#triang  &  iter  &  time  & iter  &  time  & iter  &  time \\
		   1089 &    2048  &   119  &  0.10  &   68  &  0.04  & 1000  &  0.76 \\
		   4225 &    8192  &   112  &  0.20  &   74  &  0.10  &  299  &  0.45 \\
		  16641 &   32768  &    57  &  0.41  &   80  &  0.54  &  104  &  0.66 \\
		  66049 &  131072  &    52  &  1.71  &   79  &  2.34  &   70  &  2.08 \\
		 263169 &  524288  &    44  &  7.45  &   78  & 11.89  &   64  &  9.97 \\
		1050625 & 2097152  &    39  & 39.00  &   79  & 74.43  &   48  & 45.43 \\
		\bottomrule
	\end{tabular}
	\caption{Numerical results in the linear case with Neumann boundary conditions}
	\label{tab:linear_nbc}
\end{table}
\begin{table}
	\centering
	\begin{tabular}{rrrrrrrr}
		\toprule
		& & \multicolumn{2}{c}{PG $L^1(\Omega)$} & \multicolumn{2}{c}{PG $L^2(\Omega)$} & \multicolumn{2}{c}{Cond.grad.}\\
		\cmidrule(r){3-4} \cmidrule(r){5-6} \cmidrule(r){7-8}
		\#nodes &\#triang  &  iter &   time  & iter  &  time  & iter  &  time  \\
		   1089 &    2048  &   127 &   0.69  &   77  &   0.42 &  1000 &   4.70 \\
		   4225 &    8192  &    56 &   1.40  &   59  &   1.32 &   172 &   3.51 \\
		  16641 &   32768  &    83 &   7.78  &   69  &   6.39 &    96 &   9.05 \\
		  66049 &  131072  &    35 &  19.97  &   66  &  31.24 &    72 &  29.49 \\
		 263169 &  524288  &    37 & 104.96  &   66  & 145.05 &    56 & 125.62 \\
		\bottomrule
	\end{tabular}
	\caption{Numerical results in the nonlinear case with Neumann boundary conditions}
	\label{tab:nonlinear_nbc}
\end{table}

\printbibliography
\end{document}